\documentclass{article}
\usepackage{latexsym}
\usepackage{graphics}
\usepackage{amsmath}
\usepackage{amsfonts}
\begin{document}
\bibliographystyle{plain}
\title{On conformally invariant differential operators.}
\author{Spyros Alexakis}
\maketitle
\newtheorem{proposition}{Proposition}
\newtheorem{theorem}{Theorem}
\newtheorem{lemma}{Lemma}
\newtheorem{conjecture}{Conjecture}
\newtheorem{observation}{Observation}
\newtheorem{formulation}{Formulation}
\newtheorem{definition}{Definition}
\newtheorem{corollary}{Corollary}

\begin{abstract}
We construct new families of conformally invariant differential
operators acting on densities. We introduce a simple, direct
approach which shows that all such operators arise via this
construction when the degree is bounded by the dimension. The
method relies on a study of well-known transformation laws and on
Weyl's theory regarding identities holding ``formally'' vs. ``by
substitution''. We also illustrate how this new method can
strengthen existing results in the parabolic invariant theory for
conformal geometries.
\end{abstract}

\section{Introduction}

\par This paper presents a construction of new families of conformally invariant differential
operators acting on densities and partially  shows that all such
invariant operators arise via this construction. This project thus
fits into the program of identifying invariants of parabolic
geometries, a problem on which there is a rich literature and for
which an invariant theory has been developed (see \cite{eg:icd},
\cite{beg:itccg},\cite{{g:itccg}}, \cite{{gh:cipl}} and references
therein). In this paper we restrict attention to conformal geometry
(and not CR or projective geometry).

\par The broad challenge of constructing local objects (scalars,
tensors, differential operators etc) which exhibit a form of
invariance under conformal changes of the underlying metric has
been pursued for some time, partly in connection with questions in
general relativity, see \cite{t:cg}, \cite{s:ct}, \cite{p:cti}. A number of closely
related techniques for constructing such objects have been developed, e.g. the works of T.Y.
Thomas \cite{t:digs}, the Cartan conformal connection (see \cite{k:tgdg}), the Fefferman-Graham {\it
ambient metric} \cite{fg:ci}, and the tractor calculus \cite{beg:tsbcprs}. Our construction uses
the ambient metric.

\par The method we employ for the construction of the new
operators is the standard method with which one obtains
conformally invariant scalar quantities using the ambient metric:
 Very roughly, the ambient metric provides a way to embed
an entire conformal class $(M,[g])$ (and also a conformal density
defined over it) into an ambient pseudo-Riemannian manifold, so
that any intrinsic scalar object one constructs in the ambient
manifold will automatically be a conformal invariant of the
original conformal class $(M,[g])$. Now, once one constructs
conformally invariant objects (differential operators in this case),
 it is a natural question to ask
whether one has found all such objects which exhibit the required
invariance. This can be thought of as a completeness question.

\par Our proof that  all operators with the
 required conformal invariance arise via our construction
(theorems \ref{ayto}, \ref{ayto2}, \ref{careful}) is direct and in
a sense elementary. It presupposes no knowledge of representation
theory--in particular no prior knowledge of \cite{eg:icd},
\cite{beg:itccg} is needed. Essentially, our proof relies on a
careful study of the transformation laws under conformal
re-scaling of the objects under consideration and the classical
formalism of Weyl \cite{w:cg}. Our method has one limitation: It
can only be applied if the degree $d$ (see Definition
\ref{kakom}) is less than or equal to the dimension $n$.
 For completeness, we also illustrate how the invariant theory
developed in \cite{beg:itccg} can be applied to settle the case
where the degree is greater than the dimension, provided $n$ is
odd and the densities satisfy certain additional restrictions
 (see theorem \ref{easy} below).

\par It is worth noting that a forthcoming paper by Hirachi
(see \cite{h:new}) develops an analogous
  direct argument to construct and prove completeness for CR-invariant operators, with
   applications to the asymptotic expansion of the Szeg\"o
    kernel of strictly pseudo-convex domains.  It is also interesting to note
that our approach is in some sense analogous to the ``direct
method'' that was used in Proposition 3.2 in \cite{beg:itccg} for
the case of invariants of densities with degree bounded by the
dimension. The authors in \cite{beg:itccg} used this method in the
case of invariants involving only densities,  not  curvature
tensors (the latter case is more subtle due to the algebraic
complexity of the curvature and its covariant derivatives).
\newline

{\it Outline of the paper:} In section 2 we give a rigorous
definition of conformal invariance for differential operators in
curved space, and we recall the earlier known examples of such
operators. In section 4 we recall the Fefferman-Graham ambient
metric construction and explain in detail the sense in which this
is a conformally invariant construction. We then construct the new
operators and discuss some of their features (they arise in
families and most of them vanish in flat space). In section 5 we
state the completeness results, which we then prove in sections 6
and 7. In section 6 we first give a one-page synopsis of the
argument, and then explain how the proof is naturally divided into
three steps; we prove the first two (which are shorter) in the
rest of section 6 and the lengthier one in section 7. Finally, in
section 8 we present a straightforward adaptation of the methods
in \cite{beg:itccg} to prove a completeness result in the case where
the degree is higher than the dimension (in which case our new
approach does not go through).
\newline

{\it History of the problem and strengthening of existing results:}
 The first pioneering construction of conformally invariant scalars
(depending only on curvature terms) was carried out in
\cite{fg:ci}, where Fefferman and
 Graham posed the geometric problem whether all invariants arose via their construction.
The natural notion of conformal invariance for general curved
structures is the one discussed in
Definition \ref{defconfinv} below. We will be using this notion of
invariance throughout our paper.

\par The papers \cite{eg:icd}, \cite{beg:itccg} sought to address the problem
 posed in \cite{fg:ci}.  The paper \cite{beg:itccg} of Bailey,
 Eastwood and Graham was the first to attack
 the geometric problem where (subject to certain restrictions in even dimensions)
the authors proved completeness for conformal invariants which
 locally depend on the curvature and its derivatives.  In \cite{eg:icd} and also in \cite{beg:itccg}
the authors also discuss the problem of determining all conformal
 invariants that depend on the derivatives of a single function
 (and these are called conformal invariants of densities).
This latter question is studied as an
interesting model problem, which does not apply directly to the
geometric problem, see \cite{g:itpg}.

\par Now, the problem that we address in this paper is to
construct and prove completeness for invariant {\it differential
operators} (which locally depend on {\it both} curvature terms and
derivatives of a function). Nonetheless, our approach to this
problem succeeds in settling some cases
 that \cite{eg:icd} and \cite{beg:itccg} left open:

 In even dimensions, the authors in \cite{beg:itccg}
solve the completeness problem for conformally invariant scalars
provided the degree is strictly less than $n$ ($n$ being
the dimension). Our method also captures the case where the degree
is $n$.

\par Regarding conformal invariants of densities, the authors in
\cite{eg:icd} and \cite{beg:itccg} consider the problem in the
setting of Euclidean space, where they require invariance of the
local objects under the action of the conformal group. This is different
 from our setting, where we consider objects defined on all curved
structures and the notion of conformal invariance we require is
the one of definition \ref{defconfinv}. Nonetheless, in the
setting of Riemannian operators which are conformally invariant in
the sense of definition \ref{defconfinv} (which is a stronger
requirement than the invariance under the conformal group imposed
in \cite{eg:icd}),
 we can strengthen some of the results in \cite{eg:icd}, \cite{beg:itccg} regarding
invariants of densities. See the discussion in section 5.
\newline

{\it Applications:} Fefferman and Hirachi \cite{fh:amcqcccg} have
proven that each new operator $P_{g}$ governs the transformation
law (under conformal re-scalings) of a certain scalar
$Q^{P_{g}}$ which depends on both the Weyl and the Ricci
curvatures. Thus, the study of these new scalars can be seen as
the study of the interplay between the Weyl and Ricci curvatures
in a conformal class.

\section{Formulation of the problem.}

\par Our goal here is to explicitly construct all the {\it intrinsic}
differential operators $L_{g}(f)$, defined on smooth
 Riemannian manifolds $(M,g)$ of a fixed dimension $n$, that remain
invariant under conformal re-scalings of the metric $g$. We first define
 {\it intrinsic} operators\footnote[1]{Our notion of intrinsicness in definition \ref{intrinoper} requires
 the symbol of the operator to remain invariant under both orientation-preserving
 and orientation-reversing isometries, Thus, we are actually restricting
attention to {\it even} invariants, in the language of \cite{beg:itccg}}.

\begin{definition}
\label{intrinoper}
 An {\it intrinsic} differential operator $L_{(M,g)}$ acting on scalar functions $f\in C^\infty(M)$ is a
differential operator associated to each Riemannian manifold $(M,g)$ with the following properties:

\begin{enumerate}

\item{There is a fixed polynomial expression $P$ in the variables
$\frac{\partial^a g_{ij}}{\partial x_1^{k_1}\dots
\partial x^{k_n}_n}, (det g_{ij})^{-1}$ and $\frac{\partial ^b
f}{\partial x_1^{l_1}\dots \partial x^{l_n}_n}$ so that for any
Riemannian manifold $(M,g)$ of dimension $n$ and any $f\in
C^\infty(M)$ we have in local coordinates the formula

$$L_{(M,g)}f=P(\frac{\partial^a g_{ij}}{\partial x_1^{k_1}\dots \partial x^{k_n}_n},
 (det g_{ij})^{-1} ,\frac{\partial ^b f}{\partial x_1^{l_1}\dots \partial x_n^{l_n}})$$ }

\item{For any two Riemannian manifolds $(M,g), (M',g')$, isometric
via the map $\phi:M\rightarrow M'$, and any $f\in C^\infty (M)$,
we have that: $L_{(M',g')}\big{(}\phi_{*}
f\big{)}=\phi_{*}(L_{(M,g)}f)$.}
\end{enumerate}
\end{definition}

\par It is a classical result of Weyl \cite{w:cg} (see also \cite{e:ntrm})
that any intrinsic (also called  Riemannian) differential operator
can
 be written as  a linear combination of complete contractions
 involving intrinsic ``building blocks'':

\par Let $R$ stand for the curvature tensor of the metric
$g$ (we are suppressing the lower indices of the curvature
tensor).
 Let $\nabla$ stand for the Levi-Civita connection of
$g$. Weyl's theory then ensures that for every {\it intrinsic}
operator $L_g(f)$, there exists a fixed linear combination
$\sum_{x=1}^V a_x C^x_g(f)$ of complete contractions in the form:

\begin{equation}
\label{contr1} contr(\nabla^{(m_1)}_{a_1\dots
a_{m_1}}R_{ijkl}\otimes\dots\otimes \nabla^{(m_s)}_{b_1\dots
b_{m_s}}R_{i'j'k'l'}\otimes \nabla^{(p_1)}_{h_1\dots
h_{p_1}}f\otimes\dots\otimes \nabla^{(p_q)}_{z_1\dots z_{p_q}}f)
\end{equation}
(the contractions of indices are being taken with respect to the metric $g$) so that
for every Riemannian manifold $(M,g)$ and every $f\in C^\infty
(M)$:

\begin{equation}
\label{tampela} L_g(f)= \sum_{x=1}^V a_x C^x_g(f) \end{equation}

\par In order to define {\it conformally invariant}
differential operators, we also need a notion of weight for the
operators in the  form (\ref{tampela}): For any complete
contraction in the form (\ref{contr1}), we define its {\it weight}
$K$ to be the number $K=-\sum_{k=1}^s (m_k+2)-\sum_{k=1}^q p_k$.
Observe that complete contractions $C_g(f)$ with weight $K$ have
the property that for every $t\in \mathbb{R}$ they must satisfy
$C_{t^2g}(f)=t^{K}C_{g}(f)$.
 \newline

\begin{definition}
\label{operweight} An intrinsic operator $L_g(f)$ that can be
written as a linear combination $L_g(f)=\sum_{x=1}^V a_x C^x_g(f)$
 of complete contractions in the form (\ref{contr1})
 with weight $-K$ will be called an intrinsic operator of weight $-K$.
\end{definition}

\par Now, recall that two metrics $g,g'$ defined over a manifold $M$ are
{\it conformally equivalent} if there exists a function $\phi\in
C^\infty (M)$ so that $g=e^{2\phi}g'$. Now given any metric $g$
defined over a manifold $M$, the set of metrics $g'$ that are
conformally equivalent to $g$
 define an equivalence class, which we denote by $[g]$.
\newline

 A natural question is then to determine which
Riemannian differential operators $L_g(f)$
exhibit invariance properties under conformal transformations of
the metric $g$. In the setting of differential operators, the
natural generalization of the function space $C^\infty (M)$ is the
bundle $E[w]$ of conformal densities of a given weight $w\in
\mathbb{R}$.

\begin{definition}
\label{confdens}
 Given any conformal class $(M,[g])$ we define a $w$-density
$f_w$ (of weight $w$) to be a function

$$f_w: M\times [g]\longrightarrow \mathbb{R}$$
so that for any pair $g_1,g_2\in [g]$ with $g_2= e^{2\phi}g_1$ we
have that:

$$f_w(x,g_2)= e^{w\cdot \phi} f(x,g_1)$$
We denote the bundle of densities of weight $w$ by $E[w]$. (Note
that sometimes $E[w]$
 also denotes the space of sections of this bundle).
\end{definition}

\begin{definition}
\label{defconfinv} An intrinsic differential operator $L_{g} (f)$
will be called conformally invariant of bi-degree $(a,b)$ if
 for every $\Omega>0$, $\Omega\in C^{\infty}(M^n)$ we have:

\[
L_{\Omega^2g}\big(\Omega^{a}f\big)=\Omega^{b}L_{g}(f).
\]
\end{definition}
In more formal language, we can say that an intrinsic operator
$L_g$ of bi-degree $(a,b)$ maps $E[a]$ into $E[b]$, in the sense
that if $f_w$ belongs to the bundle $E[a]$ then $L_g(f_w)$ will be
an element of the bundle $E[b]$.

\par Examples of conformally invariant operators on densities have been
 known for some time.
The most classical example is the conformal Laplacian:

$${\Delta}^c_{g}:
E[-{n\over 2}+1]\to E[-{n\over 2}-1].$$

which in dimension $n$ is given by the formula:

$$\Delta^c_{g}= [\Delta_{g} +\frac{n-2}{2(n-1)}S_{g}]$$
where $\Delta_{g}$ is the Laplace-Beltrami operator and $S_{g}$ is
the scalar curvature.

\par In 1984 Paneitz showed
that for $n=4$  one can add lower order terms to ${\Delta}^2_{g}$
and make it conformally invariant of bi-degree $(0,-4)$. Branson
\cite{b:fd} later generalized this result to arbitrary
dimensions $n\ge 4$. He showed that the following operator is
conformally invariant of bi-degree $(-\frac{n}{2}+2,
-\frac{n}{2}-2)$.

\begin{equation}
\label{paneitz} P_4^n\psi = \Delta^2_{g}\psi -div (a_n S_{g}
gd\psi +b_n Ric_{g}d\psi) +\frac{n-4}{2}Q^n_{g} \psi
\end{equation}
where $a_n=\frac{(n-2)^2+4}{2(n-1)(n-2)}$ and $b_n=-\frac{4}{n-2}$
and also $Q^n_{g}= -\frac{1}{2(n-1)}\Delta_gS_g+
\frac{n^3-4n^2+16n-16}{8(n-1)^2(n-2)^2}S_g^2 -\frac{2}{(n-2)^2}
|Ric_g|^2$. Here  $Ric_g$ is the Ricci curvature.

\par Finally, the authors in \cite{gjms:cipl} proved that
in odd dimensions, one can add lower order terms to the $k^{th}$
power of the Laplacian, $\Delta^k$,
 to obtain a conformally invariant operator $P^n_{2k}$
of bi-degree $(-\frac{n}{2}+k,-\frac{n}{2}-k)$.
$P^n_2$ is then ${\Delta}^c$ and $P^n_4$ is the Paneitz operator. For
even dimensions, the same construction goes through provided $k\le \frac{n}{2}$
(see also \cite{gh:cipl} for a non-existence theorem which shows
 that this result is sharp). These operators are now called the GJMS operators.
\newline

\par It is known (see \cite{j:cc}) that in conformally flat space, the
powers of the Laplacian
 are the only nontrivial linear conformally invariant operators.
In section 4 we provide a general construction of conformally invariant operators
in general curved spaces, and in
sections 6 and 7
 we show that under certain restrictions, all invariant operators arise via our construction.

\section{Notational Conventions}

\par  $\mathbb{Z}_{+}$
will stand for the set of {\it strictly positive} integers.
 Throughout this paper $n$ will stand for the (fixed) dimension in
which we are
 considering our operators, and $g$ will stand for a
 metric tensor of an $n$-dimensional manifold $M$.
 When we
wish to consider operators for dimensions $N\ne n$, we will be
explicitly writing out $g^N$. {\it Note:} We will
be assuming that $n\ge 3$ throughout the paper. The case $n= 3$ is
slightly different from the other cases. We will be adding
special footnotes regarding the case $n=3$ whenever necessary.

\par Throughout this paper, we will be writing out linear combinations of
complete contractions; it will be useful to impose certain conventions regarding the form of the
factors in these complete contractions.

 First, when we write $\nabla^{(m)}$,
 then $m$ will stand for the number of differentiations. If we write $\nabla^a$,
then ${}^a$ will be a raised index. Furthermore, we will usually
be writing $\nabla^{(m)}R$ for the $m^{th}$ covariant derivative
of the curvature tensor without writing out the indices of this
tensor (i.e. we will not be writing out $\nabla^{(m)}_{r_1\dots
r_m}R_{ijkl}$). However, we impose the restriction that when a
factor $\nabla^{(m)}R_{ijkl}$ appears in a complete contraction,
then the indices ${}_i,{}_j,{}_k,{}_l$ are {\it not} contracting
against each other.

\par Furthermore, for any linear
combination in the form $\sum_{h\in H} a_h C^h_{g}(f)$ and any
subset $H'\subset H$, we will call $\sum_{h\in H'} a_h C^h_{g}(f)$
a {\it sublinear combination} of $\sum_{h\in H} a_h C^h_{g}(f)$.

\par Now, an important note regarding the notion of
identities ``holding'': The operators $L_g(f)$ that we will be
considering will be functions of a metric $g$ and a function $f$.
They will be written out as linear combinations of complete
contractions. Now, when we write $\sum_{h\in H} a_h
C^h_g(f)=\sum_{p\in P} a_p C^p_g(f)$, we will mean that for any
manifold $M$, any metric $g$, any $f\in C^\infty(M)$ and any
$x_0\in M$, the two sides of this equation have the same values at
$x_0$. Thus, throughout this paper, when we prove that a sublinear
combination $\sum_{h\in H'} a_h C^h_{g}(f)$ in $L_g(f)=\sum_{h\in
H} a_h C^h_{g}(f)$ is equal to some expression $\sum_{t\in T} a_t
C^t_g(f)$, we will be free to {\it replace} the sublinear
combination $\sum_{h\in H'} a_h C^h_{g}(f)$ in $L_g(f)$ by the
linear combination $\sum_{t\in T} a_t C^t_g(f)$.

\section{The construction of the new operators.}

\par The only piece of background needed for the construction
 of the operators is the {\it ambient metric}, introduced by Fefferman and Graham in \cite{fg:ci}.
The ambient metric is a formal construction that invariantly
associates to each conformal class $(M,[g])$ an
$(n+2)$-dimensional pseudo-Riemannian manifold $(\tilde{G},
\tilde{g})$ (more precisely, a jet of an $(n+2)$-metric
$\tilde{g}$).
 It was this tool
that was used, albeit in a different manner, in \cite{gjms:cipl}
to construct the conformally invariant powers of the Laplacian.

\subsection{The Fefferman-Graham ambient metric.}
\label{amb}

\par All the material presented in this subsection comes from \cite{fg:ci} and
 \cite{fh:amcqcccg}.
 Let $(M,[g])$ be a conformal class and $g$  a representative of this class.

\par Define
$G= \mathbb{R}_{+}\times M$. Any coordinate patch
$U\subset M$ (with coordinates $x^1,\dots x^n$)
 defines a coordinate patch $\mathbb{R}_{+}\times U$ in
$G$ (with coordinates $t,x^1,\dots x^n$).
 Define a symmetric $(0,2)$-tensor
  $g^{n+1}$ on $G$ via the formula:

$$\sum_{i,j=0}^n g^{n+1}_{ij}(t, x)dx^idx^j= t^2
\sum_{i,j=1}^n g_{ij}(x)dx^idx^j$$ (hence, the $t$-direction is
null). Now let $\tilde{G}=G\times (-1,1)$, where $\tilde{G}$
 has coordinates $\{t,x^1,\dots ,x^n,\rho\}$ and also $G=\{\rho=0\}$.
 Fefferman and Graham have proven that {\it if the dimension $n$ is odd} then
there exists a metric $\tilde{g}_{ij}$ (any value $k, 1\le k\le n$ of the indices $i,j$ corresponds to
the vector $\frac{\partial}{\partial x^k}$ and the values $0$ and $n+1$ correspond to the vectors
 $\frac{\partial}{\partial t}$, $\frac{\partial}{\partial \rho}$)
 on $\tilde{G}$ with the following properties:

\begin{enumerate}
\item{$\tilde{g}(t,x,0)|_{TG}=
g^{n+1}(t,x)|_{TG}$.}

\item{For $1\le i,j\le n$, $\tilde{g}_{ij}(t,x,\rho)=t^2\tilde{g}_{ij}(1,x,\rho)$.}

\item{Off of the hypersurface $\{ \rho =0\}$, we have
$Ric(\tilde{g})(t,x,\rho) =0+O(\rho^\infty)$.}
\end{enumerate}

\par Furthermore, there exists a coordinate system $\{t,x,\rho\}$ on $\tilde{G}$ for which
 the metric
$\tilde{g}_{ij}(t,x,\rho)$ can be written in a special form:
Denoting $dx^0=dt$ and $dx^{n+1}=d\rho$, we have:

\begin{equation}
\label{analysh} \sum_{i,j=0}^{n+1}\tilde{g}_{ij}dx^idx^j=
t^2\sum_{i,j=1}^n g_{ij}(x,\rho)dx^idx^j +2tdt d\rho +2\rho dt^2
\end{equation}

\par For each ambient metric construction, where we start off with
$(M,g)$ and perform the above construction, we will call this
 coordinate system $(t,x,\rho)$ the
{\it special coordinate system} that corresponds to $(M,g)$.

\par Fefferman and Graham have then shown in \cite{fg:ci} that {\it for $n$ odd},
the Taylor expansion of the metric $\tilde{g}$ off of the hypersurface
$\{\rho=0\}$ {\it is
 uniquely determined by the above requirements}. Thus, for a given
 metric $g\in [g]$, two different ambient metric
  constructions differ by terms in $O(\rho^\infty)$.
 Moreover, the ambient
 metric is a {\it conformally invariant construction}, in the
following sense: Let $g_1,g_2\in [g]$, where $g_1=e^{2\phi}g_2$.
Let $(\tilde{G}_1,\tilde{g}_1)$, $(\tilde{G}_1,\tilde{g}_2)$ be
ambient metric constructions for $g_1,g_2$. We denote by
$(t,x^1,\dots x^n,\rho)$ the special coordinate system that
corresponds to $\tilde{g}_1$ and by $(t',x^1,\dots x^n,\rho')$ the
special coordinate system that corresponds to $\tilde{g}_1$. Then
there exists a map $\Phi:\tilde{G}_2\rightarrow \tilde{G}_1$ so
that:

\begin{enumerate}
\item{ $\Phi(1,x^1,\dots ,x^n,0)=(e^{\phi(x^1,\dots ,x^n)},
x^1,\dots x^n,0)$ and $\Phi$ maps the set $\{ \rho'=0\}$ onto the set
$\{\rho=0\}$.} \item{$\Phi$ respects the homogeneity of
$\tilde{G}_1$ and $\tilde{G}_2$,
 in the sense that if $\Phi(t',x,\rho')=(t,x,\rho)$ then $\Phi(\lambda\cdot t',x,\rho')=(\lambda\cdot t,x,\rho)$.}
\item{$\big{(}\Phi^{*}\tilde{g}_1\big{)}= \tilde{g}_2 +O(\rho^\infty)$
 (so the ambient metric constructions for $g_1,g_2$
 are isometric mod $O(\rho^\infty)$).}
\end{enumerate}

\par When $n$ is even the ambient metric construction can only be carried out to finite order.
With the notational conventions introduced above, Fefferman and Graham have shown that there exists
a metric $\tilde{g}$ on $\tilde{G}$ so that:

\begin{enumerate}
\item{$\tilde{g}(t,x,0)|_{TG}=
g^{n+1}(t,x)|_{TG}$.}

\item{For $1\le i,j\le n$, $\tilde{g}_{ij}(t,x,\rho)=t^2\tilde{g}_{ij}(1,x,\rho)$.}

\item{Off of the hypersurface $G^{n+1}=\{ \rho =0\}$, we have that:
$Ric(\tilde{g})(t,x,\rho) =0+O(\rho^{\frac{n-4}{2}})$, while
components of $Ric(\tilde{g})(t,x,\rho)$ that are tangential
to $G^{n+1}$ vanish to order $\frac{n-2}{2}$.}
\end{enumerate}
Furthermore, there exists a coordinate system $\{t,x,\rho\}$
on $\tilde{G}$ for which the metric
$\tilde{g}_{ij}(t,x,\rho)$ can be written in a special form:
Denoting $dx^0=dt$ and $dx^{n+1}=d\rho$, we have:

\begin{equation}
\label{analysh'} \sum_{i,j=0}^{n+1}\tilde{g}_{ij}dx^idx^j=
t^2\sum_{i,j=1}^n g_{ij}(x,\rho)dx^idx^j +2tdt d\rho +2\rho
dt^2+O(\rho^{\frac{n}{2}})
\end{equation}

 Fefferman and Graham have then shown in \cite{fg:ci} that {\it for $n$ even},
the Taylor expansion of the metric $\tilde{g}^{n+2}$ off of the hypersurface
$\rho=0$ {\it is
 uniquely determined by the above requirements up to order $\frac{n}{2}$}.
  Thus, for a given metric $g\in [g]$, two different ambient metric constructions
   differ by terms in $O(\rho^{\frac{n}{2}})$.

\par We also note that if we denote by $\tilde{R}_{ijkl}$ the ambient
 curvature tensor  and by $\tilde{\nabla}$ the Levi-Civita
connection of the metric $\tilde{g}$, then
$\tilde{R}_{ijkl}$ is related to the
 curvature tensor $R_{ijkl}$ of the underlying manifold
 $(M,g)$ in a simple way:

\par Recall the Schouten and Weyl tensors:

\begin{equation}
\label{schouten} P_{ij}=\frac{1}{n-2}[Ric_{ij}-\frac{S}{2(n-1)}
g_{ij}]
\end{equation}

\begin{equation}
\label{weyl}
W_{ijkl}=R_{ijkl}-[P_{jk}g_{il}+P_{il}g_{jk}- P_{ik}g_{jl}-
P_{jl}g_{ik}]
\end{equation}

\par Fefferman and Graham, \cite{fg:ci}, have shown that at each point $(t,x)$ of $G^{n+1}$,
  we have that for the vectors $X^0,X^1\dots X^n,X^\infty$ that correspond
   to the coordinates $t,x^1,\dots ,x^n,\rho$:

\begin{enumerate}

\item{$\tilde{R}_{ijk0}(t,x,0)=0$, $0\le i,j,k\le n+1$.}

\item{$\tilde{R}_{ijkl}(t,x,0)= t^2W_{ijkl}(x)$, $1\le i,j,k,l\le n$.}

\item{$\tilde{R}_{ijk\infty}(t,x,0)=t^2 C_{kij}(x)$, $1\le i,j,k$.}

\item{$\tilde{R}_{\infty ij\infty}(t,x,0)=\frac{t^2}{n-4}B_{ij}(x)$, $1\le
i,j \le n, n\ne 4$}
\end{enumerate}
where  $C_{kij}$ is the Cotton tensor, $C_{kij}=
\nabla_iP_{jk}-\nabla_jP_{ik}$  and $B_{ij}$ is the Bach
 tensor, $B_{ij}= {C_{ijk,}}^k -P^{kl}W_{kijl}$.

\par Moreover, it has been shown in \cite{g:ciplne} and that the
 Christoffel symbols
$\tilde{\Gamma}^k_{ij} (\tilde{x}_0)$ (with respect to the special
coordinate system) are related to the underlying geometry of
$(M,g)$ by simple relations:
 Let the indices $a,b,c$ take values between $1$ and $n$.
 Then at each point $\tilde{x}_0=(1,x_0,0)\in \tilde{G}$, $x_0\in M$, we have:

\begin{eqnarray}
\label{Gamma}
\tilde{\Gamma}^a_{bc}(\tilde{x}_0)=\Gamma^a_{bc}(x_0), \text{  }
\tilde{\Gamma}^\infty_{ab}(\tilde{x}_0)=g_{ab}(x_0), \text{  }
\tilde{\Gamma}^0_{ab}(\tilde{x}_0)=-P_{ab}(x_0)
\end{eqnarray}
 The rest of the
Christoffel symbols can be computed using the formula
$\frac{\partial \tilde{g}_{ij}}{\partial\rho}(1,x,0)=2P_{ij}(x)$,
see \cite{g:ciplne}, and the form (\ref{analysh}) of the ambient
metric at $(t,x,0)$.

\subsection{The construction of the New Operators.}

\par In this subsection we will construct the new
 operators, and explain some of their features. It should be noted that the method we use is
the standard way one constructs conformally invariant scalar
objects using the ambient metric (initiated in \cite{fg:ci}). We
will initially construct the new operators in the odd-dimensional
case when $(w+\frac{n}{2})\notin\mathbb{Z}_{+}$. We will then
explain how the method can be carried over to even dimensions
and/or $(w+\frac{n}{2})\in\mathbb{Z}_{+}$.
\newline

\par We start with any conformal equivalence class $(M,[g])$ and a density $f_w(x,g)$
of weight $w$ defined over $[g]$. We pick any fixed metric $g$
from the class $[g]$.
 Evaluating the density $f_w$ at
the metric $g$ defines a scalar function $f(x)=f_w(x,g)$, $f\in
C^\infty(M)$.

 We perform the ambient metric construction
$(\tilde{G}, \tilde{g})$ for $(M,g)$. In the ambient metric setting,
 the density $f_w$ can be naturally viewed as a homogeneous
function $u_w(t,x)$ defined on $G\subset \tilde{G}$ by setting:

\begin{equation}
\label{prodernamen} u_w(t,x)=t^w f(x)
 \end{equation}

\par We then seek to invariantly extend this homogeneous
 function $u_w(t,x)$ to a
 function $\tilde{u}_w(t,x,\rho)$ defined on
$(\tilde{G}, \tilde{g})$. We do this by requiring the extension
$\tilde{u}_w$ to be homogeneous and harmonic to infinite order off of
$\{\rho=0\}$:

\begin{enumerate}

\item{For the special coordinate system $(t,x,\rho)$, we require
that $\tilde{u}_w(t,x,\rho) =t^w \tilde{u}(1,x,\rho)$.}

\item{We require that $\tilde{u}_w(t,x,0)=u_w(t,x)$.}

\item{We require that:
\begin{equation}
\label{Laplace}
\Delta_{\tilde{g}}\tilde{u}_w(t,x,\rho)=O(\rho^\infty)
\end{equation}}
\end{enumerate}

\par Here $\Delta_{\tilde{g}}$ stands for the Laplace-Beltrami
 operator with respect to the ambient metric $\tilde{g}$ (actually if $g$
  has Riemannian signature $\Delta_{\tilde{g}}$ is the wave operator).

\par It is then known from \cite{gjms:cipl} that if
$(w+\frac{n}{2})\notin \mathbb{Z}_{+}$ then the above equation has a unique solution
$\tilde{u}_w(t,x,\rho)$, up to functions that vanish to infinite
order off of $G$. Hence, we have that the covariant
derivatives $\tilde{\nabla}^{(p)} \tilde{u}_w(t,x,0)$ at any point
$(t,x,0)\in \tilde{G}$ are all well-defined. We will refer
to the function $\tilde{u}_w$ as the {\it harmonic extension of
$u_w$ to $\tilde{G}$}.
\newline

\par Now, choose  natural numbers
$r,K\in \mathbb{N}$, and consider any complete contraction:

\begin{equation}
\label{contraction}
  contr(\tilde{\nabla}^{(k_1)}\tilde{R}\otimes
\ldots\otimes\tilde{\nabla}^{(k_s)}\tilde{R}
\otimes\tilde{\nabla}^{(l_1)}\tilde{u}_w
\otimes\ldots\otimes\tilde{\nabla}^{(l_r)}\tilde{u}_w),
\end{equation}
subject to the only restriction that $\sum l_i+\sum(k_i+2)=K$.

\par Then, for any finite set of such complete contractions,
$\{\tilde{C}^1_{\tilde{g}}(\tilde{u}_w), \dots ,
\tilde{C}^z_{\tilde{g}}(\tilde{u}_w)\}$, we can
form linear combinations in the form:

\begin{equation}
\label{lincomb2} F_{g}(f)= \sum_{s=1}^z a_s
\tilde{C}^s_{\tilde{g}}(\tilde{u}_w)
\end{equation}

\par Observe that the right hand side of the above is indeed
a function of $g$ and $f$, since for any point $x\in M$ the jets
of $\tilde{g}$ and $\tilde{u}_w$ at $(1,x,0)\in \tilde{G}^{n+2}$
are uniquely determined by the jets of $g$ and $f$ at $x\in M$.
Moreover,

\begin{proposition}
\label{propo} Let $n$ be odd and $(w+\frac{n}{2})\notin
\mathbb{Z}_{+}$. Then for any linear combination $F_{g}(f)$ as
above, $F_{g}(f)$ is an intrinsic operator of weight $-K$, which
is conformally invariant
  of bi-degree $(w, r\cdot w-K)$.
\end{proposition}

{\it Proof:} We check that $F_{g}(f)$ is an intrinsic differential
operator on $(M,g)$ by virtue of the form of the ambient curvature
tensor and the Christoffel symbols of the ambient metric (in the
special coordinate system that corresponds to $g$).

The conformal invariance  of the operator $F_{g}(f)$ follows from
the conformal invariance of the ambient metric construction, which
we discussed above: We have to show that if we pick a metric
$e^{2\phi}g\in [g]$, (for which the corresponding value of the density $f_w$ will be
$\overline{f}(x)=f_w(x,e^{2\phi}g)=e^{w\phi}f$), and we perform
the same ambient metric construction as above, then the operator
$F_{e^{2\phi}g}(\overline{f})$ will satisfy:

$$F_{e^{2\phi}g}(\overline{f})=e^{(r\cdot w-K)\cdot\phi}F_g(f). $$
In order to see this, we denote by $(\tilde{G}_2,\tilde{g}_2)$ the
ambient metric construction that corresponds to the metric
$e^{2\phi}g$, and by $\overline{u}_w$ the homogeneous and harmonic extension of
$\overline{f}$ to $\tilde{G}_2$ (as in (\ref{prodernamen})). The
discussion from the previous subsection shows us that there exists
an isometry
 $\Phi:\tilde{G}_2\rightarrow \tilde{G}_1$ with the properties
 listed in the previous subsection.
Therefore, if we denote by
$(\tilde{\nabla}^{(m)}\tilde{R})_2(1,x,0)$ the iterated covariant
derivative
 of the curvature tensor of $\tilde{g}_2$ at $(1,x,0)\in \tilde{G}_2$,
and by $(\tilde{\nabla}^{(m)}\tilde{R})_1(e^{\phi(x)},x,0))$ the
iterated covariant derivative of the curvature tensor  of
$\tilde{g}_1$ at $(e^{\phi(x)},x,0)\in \tilde{G}_1$, we will have
that:

$$(\tilde{\nabla}^{(m)}\tilde{R})_2(1,x,0)=\big{(}\Phi_{*}
(\tilde{\nabla}^{(m)}\tilde{R})_2\big{)}(e^{\phi(x)},x,0)=
(\tilde{\nabla}^{(m)}\tilde{R})_1(e^{\phi(x)},x,0)=
e^{2\phi(x)}(\tilde{\nabla}^{(m)}\tilde{R})_1(1,x,0)$$

We simlarly observe that $(\Phi_{*}\overline{u}_w)$ must still be
a homogeneous harmonic function in $\tilde{G}_1$. Therefore we
derive $\Phi_{*}\overline{u}_w=\tilde{u}_w+O(\rho^\infty)$.

\par Therefore, we have that:

$$ F_{e^{2\phi}g}(\overline{f})(x)=\sum_{s=1}^z a_s
\tilde{C}^s_{\tilde{g}_1}(\tilde{u}_w)(e^{w\phi(x)},x,0)$$
No, using the fact that $\tilde{g}^{ij}_1(t,x,0)=t^{-2}\tilde{g}^{ij}_1(1,x,0)$,
$(\tilde{\nabla}^{(m)}\tilde{R})_1(t,x,0)=t^2(\tilde{\nabla}^{(m)}\tilde{R})_1(1,x,0)$,
 $(\tilde{\nabla}^{(p)}\tilde{u}_w)_1(t,x,0)=t^w(\tilde{\nabla}^{(p)}\tilde{u}_w)_1(1,x,0)$,
 we derive our proposition. $\Box$

\begin{definition}
\label{Weylop} A differential operator in the form
(\ref{lincomb2}) will be called a Weyl operator of weight $-K$ and
$f$-homogeneity $r$ in $\tilde{u}_w$.
\end{definition}

{\it The case of half-integer weights and/or even dimensions:}
 If $n$ is odd, then for each weight $w=-\frac{n}{2}+k,
k\in \mathbb{Z}_{+}$ it is shown in \cite{gjms:cipl} that the
equation (\ref{Laplace}) has a uniquely defined
 solution up to the order $k-1$. Thus the same proof as above shows that
the operators in (\ref{lincomb2})
 are well-defined and conformally invariant, provided that for each
$\tilde{C}^s_{\tilde{g}}(\tilde{u}_w)$ (in the form
(\ref{Wcontr})) we have that $l_h\le k-1$.

\par If $n$ is even, then we have seen that the ambient metric is
 well-defined only up to order $\frac{n}{2}$. Hence, the same proof as above shows that for
 each weight $w, (w+\frac{n}{2})\notin\mathbb{Z}_{+}$, the operators in (\ref{lincomb2}) are well-defined
 provided that, when the expressions (\ref{contraction}) are written
 out in terms of coordinate derivatives of the ambient metric,  no derivatives
$\frac{\partial^{a_1}\dots\partial^{a_n}\partial^J
\tilde{g}^{n+2}}{\partial x^{a_1}_1\dots \partial
x_n^{a_n}\partial\rho^J}$ with $J\ge \frac{n}{2}$ appear. In the
case where $n$ is even and the weight $w=-\frac{n}{2}+k$, $k\in
\mathbb{Z}_{+}$, it follows that the operators in (\ref{lincomb2})
are well-defined provided the above restriction holds
 and also provided that for each
$\tilde{C}^s_{\tilde{g}}(\tilde{u}_w)$ we have
$l_h\le k-1$.
\newline

\par In section \ref{secthm} we will claim that all conformally invariant
differential operators arise via the above
 construction, subject to some restrictions that we will explain.
  For now, we will discuss certain interesting features of the new operators:

\subsection{Features and Examples.}

\par Let us first illustrate how each expression (when it is well-defined):

$$F_{g}(f)= \sum_{s=1}^z a_s
\tilde{C}^s_{\tilde{g}}(\tilde{u}_w)$$ can be seen as a
1-parameter family of operators, parametrized by the weight $w$.
We pick a fixed metric $g\in [g]$, and for each $w,
(w+\frac{n}{2})\notin\mathbb{Z}_{+}$ we construct the operator
$F_{g}(f)$ above (it will be conformally invariant of bi-degree
 $(w,r\cdot w-K)$.  We claim that
$F_{g}(f)$ can be expressed in the form:

\begin{equation}
\label{family} \sum_{s=1}^z a_s
\tilde{C}^s_{\tilde{g}}(\tilde{u}_w)= \sum_{h\in H} b_h(w,n)
C^h_{g}(f)
\end{equation}
where each $C^h_{g}(f)$ is Riemannian operator of weight $-K$ and
$f$-homogeneity $r$, while $b_h(w,n)$ is a rational function in
the weight $w$ and the dimension $n$.

\par {\it Proof:} This follows by calculating the Taylor expansion of $\tilde{u}_w$ off of $\{\rho=0\}$.
As discussed in \cite{gjms:cipl}, it follows that the different components
$\tilde{\nabla}^{(m)}\tilde{u}_w$ have coefficients that are rational functions in $n,w$.
\newline

\par We furthermore notice that each complete contraction
in the form (\ref{contraction}) with $s>0$ will vanish if $g$ is
locally conformally flat, since in that case
 $\tilde{R}_{ijkl}=0$ (see \cite{fg:ci}). The only complete contractions of the
form (\ref{contraction}) that do not vanish in conformally flat
space are
 the ones for which $s=0, p>1$. The operators that arise thus are
 nonlinear and have already appeared in \cite{b:tdcpitci}.

\par Moreover, since $\tilde{u}_w$ is defined to be
harmonic and $\tilde{R}$ Ricci-flat, not all the contractions
$(\ref{contraction})$ are non-zero.
 On the other hand, (for $n$ odd and $(w+\frac{n}{2})\notin \mathbb{Z}_{+}$, for simplicity)
 we can easily construct a
{\it non-zero} operator  with a leading order symbol
$C^h(\tilde{g}) \Delta^r_{g}f$  where $C^h(\tilde{g})$
is one of the conformally invariant scalars originally constructed
in \cite{fg:ci}:

\par Consider any linear combination $P(g)=\sum_{s=1}^z a_s \tilde{C}^s(\tilde{g})$
 where each $\tilde{C}^s(\tilde{g})$ is in the form:
$$contr(\tilde{\nabla}^{(m_1)}\tilde{R}\otimes\dots \tilde{\nabla}^{(m_s)}\tilde{R})$$
 with a given weight $-K$.
 These Riemannian scalars are the original conformal invariants
constructed by Fefferman and Graham in \cite{fg:ci}. Now,
for each contraction $\tilde{C}^s(\tilde{g})$, we let
$\tilde{C}^s_{\tilde{g}}(\tilde{u}_w)$ stand for the complete contraction:
$$contr((\tilde{\nabla})^{t_1\dots t_r}
[\tilde{C}^s(\tilde{g})]\otimes
 (\nabla)^{(r)}_{t_1\dots t_r}\tilde{u}_w)$$
By \cite{gjms:cipl} it follows that $\tilde{\nabla}^{(r)}_{\infty\dots \infty}\tilde{u}_w=C_{w,n}\Delta^r_gf+(lot)$,
here $C_{w,n}$ is a non-zero constant. Also, we have that each other component of
$\tilde{\nabla}^{(r)}\tilde{u}_w$ has order less than $2r$. Lastly, using the fact that
$\tilde{\Gamma}_{00}^K=0$ (and this follows from (\ref{analysh})), we see that
$\tilde{\nabla}^{(r)}_{00\dots 0}[\tilde{C}^s_{\tilde{g}}]=(-K)^r\cdot [\tilde{C}^s_{\tilde{g}}]$.
 Thus we derive our claim.
\newline

\par To illustrate, we write out two examples of new conformally
 invariant differential operators arising from the formula (\ref{contraction}).
 For both these examples we refer the reader to  \cite{fh:amcqcccg} for explicit computations.

Firstly, consider the ambient complete contraction:

$$L^1_{g}(f)=contr(\tilde{\nabla}^{m}|\tilde{R}_{ijkl}|^2\otimes
\tilde{\nabla}_m\tilde{u}_0)$$

Fefferman and Hirachi calculate that $L^1_{g}(f)$ is of the form:

$$L^1_{g}(f)= \nabla^i(|W|^2)\nabla_i f +\frac{4}{n-2}|W|^2
\Delta f$$

\par Our second example will illustrate that not all our new
operators have a leading order term $C^h(\tilde{g})
\Delta^r_{g}f$. We consider the operator $L^\sharp_{g}(f)$ that
arises by the complete contraction:

$$contr(\tilde{\nabla}^2_{st}\tilde{u}_0\otimes {\tilde{R}^s}_{jkl}
\otimes\tilde{R}^{tjkl})$$

Fefferman and Hirachi show that $L^\sharp_{g}(f)$ can be written out explicitly in the form:

\begin{equation}
\label{opallo}
L^\sharp_{g}(f)={W_{ijk}}^lW^{ijkm}\nabla^{(2)}_{lm} f-
2C_{kij}W^{ijkl}\nabla_l f+\frac{1}{n-2}|W|^2 \Delta f
\end{equation}
and hence its leading order term is
$${W_{ijk}}^lW^{ijkm}\nabla^{(2)}_{lm}  f+\frac{1}{n-2}|W|^2
\Delta f.$$

\section{The completeness Results.}
\label{secthm}

\par In this section we will present our claims that, under
certain restrictions, all conformally invariant operators can be
written as Weyl operators.
 The main restriction we need is that the weight (or the
degree for theorem \ref{careful}) be bounded by the dimension. In
the case of even dimensions and/or half-integer weight, we will
also impose additional restrictions. Since these additional
restrictions are quite technical (they depend on the parameters
$\iota,\tau, \beta,\gamma$ introduced below), the reader may wish
to skip the discussion of the extra restrictions at first.

 All three completeness
theorems \ref{ayto}, \ref{ayto2} and \ref{careful} are proven by a
novel approach that we introduce in sections 6, 7. Theorem
\ref{easy} below deals with the case where the weight is greater
than the dimension; it will be proven in section 8 by a
straightforward adaptation of the methods in \cite{beg:itccg}.
\newline

\par In order to make our task easier, we
will make an observation: Let us suppose that $L_{g}(f)=\sum_{h\in
H} a_h C^h_g(f)$ is confromally invariant of bi-degree $(a,b)$.
Let us break the index set $H$ into subsets $H_z$ according to the
rule: $h\in H_z$ if and only if $C^h_{g}(f)$ has $q=z$ (i.e. is
homogeneous of degree $z$ in the function $f$). Accordingly, we
define:

\begin{equation}
\label{fz} L^z_{g}(f)=\sum_{h\in H_z} a_h C^h_{g}(f)
\end{equation}

\par Just by applying the definition of conformal invariance
we see that:

\begin{lemma}
\label{simplicius}
\par In the above notation, if $L_{g}(f)$ is conformally invariant
of bi-degree $(a,b)$, then each $L^z_{g}(f)$ is conformally
invariant of bi-degree $(a,b)$.
\end{lemma}

\par The above Lemma allows us to restrict our attention to linear
combinations:

$$L_{g}(f)=\sum_{h\in H} a_h C^h_{g}(f)$$
where each $C^h_{g}(f)$ is in the form (\ref{contr1}) and has a
fixed homogeneity $\kappa$ in $f$. Thus, from now on we will
assume that each $C^h_{g}(f)$ has $\kappa$ factors in the form
$\nabla^{(p)}f$.

\begin{definition}
\label{kakom}
 For each complete contraction $C_{g}(f)$ in the form (\ref{contr1}) we
define $\kappa^\sharp$ to stand for the number of factors
$\nabla^{(p)}f$ with $p\ge 1$. Recall $s$ stands for the number of
factors $\nabla^{(m)}R$. We then define $2s+\kappa^\sharp$ to be
the degree of $C_{g}(f)$, $deg[C_{g}(f)]$.
\end{definition}

\par It is important to observe that any complete contraction in
the form (\ref{contr1}) with weight $-K$ will satisfy
$deg[C_{g}(f)]\le K$. Notice also that this definition is slightly different
 from the one in \cite{beg:itccg}.

\par For any complete contraction we will be paying
 special attention to pairs of indices that belong to the same
 factor and are contracting against each other. We call such
pairs of indices {\it internal contractions} (they are called
internal traces in \cite{beg:itccg}). Also, for any complete
contraction $C_g(f)$ in
 the form (\ref{contr1}) and any factor $F$ in $C_g(f)$, $\tau[F]$ will stand for the total number of
internal contractions in $F$ and $\iota[F]$ will stand for the
total number of indices in $F$ (also counting the pairs of indices
that are involved in internal contractions).

\begin{definition}
\label{betagamma} Consider any complete contraction $C_{g}(f)$ in
the form (\ref{contr1}) and we pick out its $\kappa$ factors
$F_1,\dots ,F_\kappa$ in the form $\nabla^{(p)}f$. For each
$F_s,1\le s\le\kappa$ we define $\beta[F_s]=\iota[F_s]-\tau[F_s]$.
We then define $\beta[C_{g}(f)]$ to stand for the maximum among
the numbers $\beta[F_1],\dots ,\beta[F_\kappa]$. Finally, for a
linear combination $L_{g}(f)=\sum_{h\in H} a_h C^h_{g}(f)$ we
define $\beta[L_{g}(f)]$ to be $max_{h\in H} \beta[C^h_{g}(f)]$.
\end{definition}

\par Our Theorem for odd dimensions is then the following:

\begin{theorem}
\label{ayto}
\par Let the dimension $n$ be odd.
 We pick any numbers $K\in 2\mathbb{Z}_{+}$ with $K\le n$,
$\kappa\in \mathbb{Z}_{+}$ and any weight $w$. If
$(w+\frac{n}{2})\notin \mathbb{Z}_{+}$ then any Riemannian
differential operator $L_{g}(f)$ of weight $-K$, $f$-homogeneity
$\kappa$ which is conformally invariant of bi-degree $(w, w\cdot
\kappa-K)$ can be written as a Weyl operator.

\par If $w=-\frac{n}{2}+k$ for some $k\in \mathbb{Z}_{+}$
then our conclusion above holds under the extra assumption that
$\beta[L_{g}(f)]<k$.
\end{theorem}

{\it Even dimensions:} In order to state our theorem in this case, we
introduce one more notational convention.

\begin{definition}
\label{bonjour} Given any complete contraction $C_{g}(f)$ in the
form (\ref{contr1}), we list all its factors $F_1,\dots ,F_d$. For
a factor $F_s$ in the form $\nabla^{(m)}R_{ijkl}$, we define
$\gamma[F_s]=\iota[F_s]-\tau[F_s]-2$.
 For a factor $F_s$ in the form
$\nabla^{(p)}f$, we define $\gamma[F_s]=\beta[F_s]$.

Then, for each complete contraction $C^h_{g}(f)$ we define
$\gamma[C^h_{g}(f)]$ to stand for the maximum among the numbers
$\gamma[F_1],\dots ,\gamma[F_d]$. If $L_{g}(f)=\sum_{h\in H} a_h
C^h_{g}(f)$, we define $\gamma[L_{g}(f)]=max_{h\in H}
\gamma[C^h_{g}(f)]$.
\end{definition}

\begin{theorem}
\label{ayto2}
Let the dimension $n$ be even.
 We pick any numbers $K\in 2\mathbb{Z}_{+}$ with $K\le n$,
$\kappa\in \mathbb{Z}_{+}$ and any weight $w$. If
$(w+\frac{n}{2})\notin \mathbb{Z}_{+}$,
  then any Riemannian differential operator
$L_{g}(f)$  of weight $-K$ and $f$-homogeneity $\kappa$, with
$\gamma[L_{g}(f)]< \frac{n}{2}$ which is conformally invariant of
bi-degree $(w, \kappa\cdot w -K)$
 can be written as a Weyl operator.

\par In the case where $w=-\frac{n}{2}+k$,
$k\in\mathbb{Z}_{+}$, we have that the above conclusion holds
under the additional assumption that $\beta[L_{g}(f)]< k$.
\end{theorem}

\par Let us observe that an easy consequence of the above is that
any linear conformally invariant operator $L_{g}(f)$ of bi-degree
$(-\frac{n}{2}+k,-\frac{n}{2}-k)$ with $k\le \frac{n}{2}$ can be
written in the form:

$$L_{g}(f)=(Const)\cdot P^n_{g}(f)+\sum_{h\in H} a_h
 \tilde{C}^h_{\tilde{g}^{n+2}}
(\tilde{u}_{-\frac{n}{2}+k})$$ where $P^n_{g}(f)$ is the GJMS
operator constructed in \cite{gjms:cipl}, with leading symbol
$\Delta^k$ and $\sum_{h\in H}\dots$ is a linear Weyl operator of
bidegree $(-\frac{n}{2}+k,-\frac{n}{2}-k)$: We only have to
observe that all complete contractions $C_g(f)$ in the form
(\ref{contr1}) (linear in $f$) with weight $-2k(\ge -n)$ and $s\ge
1$ factors $\nabla^{(m)}W$ will automatically have
$\gamma[C_g(f)],\beta[C_g(f)]< k\le\frac{n}{2}$.
\newline

\par Both our Theorems above require that $K\le n$, where $K$ is the weight
 of the operators and $n$ the dimension. If we wish to overcome this
restriction, however, we can prove a weaker result, which will show
 that conformally invariant operators in general can be written as Weyl
 operators, plus corrections with degree that is greater than the dimension:

 Consider any
 Riemannian operator $L_{g}(f)=\sum_{h\in H} a_h C^h_{g}(f)$,
where each $C^h_{g}(f)$ is in the form (\ref{contr1}) with weight
$-K$ and $f$-homogeneity $\kappa$. To simplify our claim,
 we will assume that each factor $\nabla^{(p)}f$ has $p>0$
  (although this restriction can easily be overcome).
For each such $L_{g}(f)$, we define its {\it minimum degree} to be
$min_{h\in H} deg[C^h_{g}(f)]$ (see definition \ref{betagamma}),
which we denote by $mindeg[L_{g}(f)]$.
 Notice that if the minimum number of factors among the contractions
 $\{C^h_{g}(f)\}_{h\in H}$ is $\sigma$, then $2\sigma+\kappa=mindeg[L_{g}(f)]$
 (i.e. the minimum degree is essentially
 determined by the minimum number of factors).
 We can then show the following:

\begin{theorem}
\label{careful} Consider a Riemannian differential operator
$L_{g}(f)$ of weight $-K$, $K\in 2\mathbb{Z}_{+}$ and
$f$-homogeneity $\kappa$. Assume $L_{g}(f)$ is conformally
invariant of bi-degree $(w,\kappa\cdot w-K)$,
$(w+\frac{n}{2})\notin\mathbb{Z}_{+}$.

\par Assume that $mindeg[L_{g}(f)]\le n$.
 Our conclusion is then that if $n$ is odd, $L_{g}(f)$ can be written as a
Weyl operator, modulo correction terms with more factors (and therefore higher degree):
 We denote by $\tilde{u}_w$ the harmonic extension of the density $f_w$.
  Then we can write:

\begin{equation}
\label{symperasma} L_{g}(f)= \sum_{h\in H'} a_h
\tilde{C}^h_{\tilde{g}^{n+2}} (\tilde{u}_w)+ \sum_{t\in T} a_t
C^t_{g}(f)
\end{equation}
where each $\tilde{C}^h_{\tilde{g}}(\tilde{u}_a)$ is a Weyl
operator and each $C^t_{g}$ has degree$>mindeg[L_{g}(f)]$.

\par In the case where $n$ is even and/or where $w=-\frac{n}{2}+k$ for
some $k\in \mathbb{Z}_{+}$, the above is still true provided
$\gamma[L_{g}(f)]< \frac{n}{2}$ and/or $\beta[L_{g}(f)]< k$.
\end{theorem}

\par We observe that this third theorem applies even when
$K>n$, provided that $mindeg[L_{g}(f)]\le n$. Thus, the theorem
 above can be iteratively applied, until we reach some linear
 combination $\sum_{t\in T'} a_t C^t_{g}(f)$ on the right
  hand side of (\ref{symperasma}) for which
  $mindeg[\sum_{t\in T'} a_t C^t_{g}(f)]>n$.
\newline

\par For all three theorems above, we will refer to the
restrictions on $\beta[C^h_{g}(f)]$ and $\gamma[C^h_{g}(f)]$
(whenever we do impose restrictions on these parameters) as
 the {\it extra restrictions}. Also, we note that in the case $n=3$,
  theorems \ref{ayto}, \ref{ayto2} can only be applied in the
 case where $L_g(f)=a_1C^1_g(f)+a_2C^2_g(f)$,where
 $C^1_g(f)=|\nabla f|^2$ and $C^2_g(f)=\Delta f\cdot f$ (because
  of the weight restrictions). On the
 other hand, still in the case $n=3$, theorem \ref{careful}
 only applies when $\kappa\le 3$ (because of the degree restriction). In the cases
 $\kappa=2,\kappa=3$, $L_g(f)$ can only be a linear combination
 without curvature terms (since $mindeg[L_g(f)]\le 3$). In the
 case $\kappa=1$, then theorem \ref{careful} only applies for
 $mindeg[L_g(f)]=3$, in which case the terms of minimum degree in
 $L_g(f)$ must be in the form
 $contr(\nabla^{(m)}R\otimes\nabla^{(p)}f)$. The cases $n=3$,
 $\kappa\ge 2$ will follow by the general argument to follow. The
 subcase $\kappa=1$ is slightly different, and we highlight this
 in footnotes, whenever needed.
\newline

{\it Theorem \ref{careful} and \cite{eg:icd}, \cite{beg:itccg}:}
A differential operator $L(f)$ (of $f$-homogeneity $\kappa$ and weight $-K$),
defined only on Euclidean space and assumed
 to be invariant under the action of the orthogonal
 group can be written as a linear combination of contractions in the form:
\begin{equation}
\label{eagles} contr(\nabla^{(p_1)}f\otimes\dots\otimes
\nabla^{(p_\kappa)}f)
\end{equation}
The additional invariance required of $L(f)$ in \cite{eg:icd} and
\cite{beg:itccg}
 is invariance under the conformal group. For comparison, we will consider
a Riemannian differential operator $L_g(f)$, whose {\it principal symbol} agrees with that of $L(f)$.
However, $L_g(f)$ may also contain complete contractions with more than $\kappa$ factors,
and these will include curvature terms. In our setting, $L_g(f)$ is required to be
conformally invariant in the sense of definition \ref{defconfinv}.

\par We observe that if we apply theorem \ref{careful}, then in Euclidean space
the correction terms will vanish, since the will involve a factor
of the curvature. Thus for $g_{Eucl}$ being the Euclidean metric,
whenever theorem \ref{careful} can be applied it will show that
$L_{g_{Eucl}}(f)$ can be written as a Weyl operator,
 without corrections terms.

\par Let us observe how Theorem \ref{careful} strengthens the existing
results for Euclidean space (in the {\it more restricted setting} of
Riemannian operators assumed to be conformally invariant in the sense of
 definition \ref{defconfinv}), provided $\kappa\le n$:
The case where $(w+\frac{n}{2})\notin \mathbb{Z}_{+}$ and $w\notin
\mathbb{Z}_{+}$ has been settled in \cite{eg:icd}. When $n$ is odd
and $w=-\frac{n}{2}+k$, theorem \ref{careful} is the first
completeness result, even in Euclidean space. In the case where
$n$ is odd and $w\in \mathbb{Z}_{+}$, \cite{beg:itccg} proves that
every invariant operator is a Weyl operator, but only provided
each $p_i\ge w$ for every $p_i$ in (\ref{eagles}). Our Theorem
\ref{careful} shows that every invariant operator is a Weyl
operator, i.e. it imposes no restriction on any $p_i$ and it works
for any $w\in \mathbb{Z}_{+}$.

For $n$ even, if $w=-\frac{n}{2}+k, k\in \mathbb{Z}_{+}$
 (which now includes the case
$w\in \mathbb{Z}_{+}$), \cite{beg:itccg} proves all invariant
operators are Weyl provided they can be expressed as polynomials
in the derivatives $\tilde{\nabla}^{(p_i)}\tilde{u}_w$, where
$\tilde{u}_w$ solves (\ref{Laplace}) and provided each $p_i\ge w$.
 Theorem \ref{careful} works without these apriori restrictions,
 and only requires $\beta[L_g(f)]<k$.

\par Finally, let us also note that our methods can also strengthen
the existing completeness results in \cite{beg:itccg} for invariants
 that depend only on the curvature (i.e. we have a linear
combination of contractions in the from (\ref{contr1}) with
$q=0$). The authors in \cite{beg:itccg} show that any such
conformal invariant must be Weyl provided the degree of the
contractions is $<n$ (in the notation of definition
\ref{kakom}--this corresponds to degree $<\frac{n}{2}$ in the
notation of \cite{beg:itccg}). Theorem \ref{careful} also shows
this claim for the case of degree $n$.
\newline

\par For completeness, we state a last theorem for the case $mindeg[L_{g}(f)]>n$,
 which will be proven in the last section, by
 an adaptation of the methods in \cite{beg:itccg}:

 \begin{theorem}
 \label{easy}
 Consider any Riemannian operator $L_{g}(f)$ with $f$-homogeneity $\kappa$, conformally
 invariant of bi-degree $(w, \kappa\cdot w-K)$.
 Suppose that $n$ is odd and $(w+\frac{n}{2})\notin
 \mathbb{Z}_{+}$,$w\notin\mathbb{Z}_{+}$, $(2\kappa\cdot
 w+n-2K)\notin-2\mathbb{Z}_{+}$, $(\kappa\cdot
 w+n-K)\notin-\mathbb{Z}_{+}$.

\par Then $L_{g}(f)$ can be written as a Weyl operator.
  \end{theorem}

\section{Proof of theorems \ref{ayto}, \ref{ayto2} and \ref{careful}: First half.}

{\it Outline of the ideas:} We prove theorems \ref{ayto},
\ref{ayto2}, \ref{careful} in steps. The main argument is
inductive: We consider the complete contractions with the smallest
number, $\sigma$,  of factors in $L_g(f)$ (denote this sublinear
combination by $L^\sigma_g(f)$) and we show that we can subtract a
Weyl operator from $L_g(f)$ to cancel out $L^\sigma_g(f)$, modulo
introducing contractions with more than $\sigma$ factors.
Iteratively repeating this step eventually shows our theorems,
since for a given weight $-K$, we cannot have more than $K$
factors in our contractions.

In order to determine the Weyl operator referred to above, and
also to prove the cancellation,
 we re-express the conformally invariant operator $L_g(f)$ as a linear
 combination of contractions involving the Weyl and symmetrized Schouten tensors (see (\ref{Wcontr})).
 It turns out that this decomposition is well-suited for our purposes, since this
naturally decomposes the curvature into a conformally invariant
and a non-conformally invariant part.

The three steps then proceed as follows: Initially we consider
 the contractions in $L^\sigma_g(f)$ and among those
 we pick out the ones with the maximum number $M$ of (symmetrized) Schouten tensors.
If $M>0$, we prove that this linear combination must vanish modulo introducing
 correction terms with $\sigma+1$ factors. This is done in Proposition \ref{queen}. By repeating
this argument, we are reduced to showing our claim when all the
contractions in $L^\sigma_g(f)$
 have no Schouten tensors. Then the key is to look at the
number of internal contractions in the terms in
$L^\sigma_g(f)$. We first pick out the complete contractions in
$L^\sigma_g(f)$ without internal contractions. We show that we can
subtract a Weyl operator from that linear combination, modulo
 introducing contractions that either have at least $\sigma+1$ factors,
{\it or} have $\sigma$ factors, no Schouten tensors, but at least
one internal contraction in some factor (Proposition
\ref{queen2}). The hardest part of the proof is to then show that
if all contractions in $L^\sigma_g(f)$ have no Schouten terms,
{\it and} have at least one
 internal contraction, then $L^\sigma_g(f)$ must vanish, modulo
correction terms with at least $\sigma+1$ factors (Proposition
\ref{queen3}). It is at this stage in the argument that we need
the degree of these contractions to be bounded by the dimension:

\par The proof of Proposition \ref{queen3} relies on
 another induction in the minimum number, $\mu$, of internal contractions
among the terms in $L^\sigma_g(f)$. We show that the sublinear combination
in $L^\sigma_g(f)$ of the
 terms with $\mu$ internal contractions
can be written as a linear combination of contractions with
 $\sigma$ factors (no Schouten tensors)
 and at least $\mu+1$ internal contractions (modulo corrections with $\sigma+1$ factors).
  In order
to prove this claim, we examine the transformation law of the
operator under conformal rescaling,
 and pick out  a very special linear combination in that transformation law that exactly corresponds to
the terms in $L^\sigma_g(f)$ with $\mu$ internal contractions.
Using Weyl's theory we deduce that this linear combination must
vanish separately.
 Then, using Weyl's theory again and an operation called ``Weylification'', we deduce that we can
make the special linear combination in $L_g(f)$ vanish, modulo
introducing terms with $\sigma$ factors and $\mu+1$ internal
contractions (and also terms with $\sigma+1$ factors). Inductively
repeating this argument we prove Proposition \ref{queen3}. (A more
detailed synopsis of this last step is provided in section 7).

\subsection{Normalizations, and the three parts of the proof.}

\par Starting with $L_g(f)$ we first re-write
$L_g(f)$ as a linear combination of contractions in a new form.

 Let us recall a
few formulas. Firstly, the curvature identity:

\begin{equation}
\label{curvature}
[\nabla_i\nabla_j-\nabla_j\nabla_i]X_l=R_{ijkl}X^k
\end{equation}

\par Secondly, we recall the Weyl and Schouten tensors (see
 (\ref{schouten}), (\ref{weyl}) above).

\par  The Weyl tensor is conformally invariant and trace-free:

\begin{equation}
\label{weyltrans} W_{ijkl}(e^{2\phi(x)}g)=e^{2\phi(x)} W_{ijkl}(g)
\end{equation}

The Schouten tensor has the following transformation law:
\begin{equation}
\label{wschtrans} P_{ij}(e^{2\phi} g)=P_{ij}(g)
-\nabla^{(2)}_{ij}{\phi}+ \nabla_i{\phi} \nabla_j{\phi}
-\frac{1}{2}\nabla^k{\phi}\nabla_k{\phi}g_{ij}
\end{equation}
while the Levi-Civita connection transforms:

\begin{equation}
\label{levicivita} ({\nabla}_k {\eta}_l)(e^{2\phi}g)=
({\nabla}_k{\eta}_l)(g) -\nabla_k\phi {\eta}_l - \nabla{\phi}_l
{\eta}_k +\nabla^s{\phi} {\eta}_s g_{kl}
\end{equation}
and the full curvature tensor $R_{ijkl}$ transforms:

\begin{equation}
\label{curvtrans}
\begin{split}
&R_{ijkl}(e^{2\phi(x)}g)=e^{2\phi (x)}[R_{ijkl}(g)+
\nabla^{(2)}_{il}{\phi}
g_{jk}+\nabla^{(2)}_{jk}{\phi}g_{il}-\nabla^{(2)}_{ik}{\phi}g_{jl}-\nabla_{jl}^{(2)}{\phi}g_{ik}
\\&+\nabla_i {\phi}\nabla_k {\phi}g_{jl}+\nabla_j {\phi}\nabla_l {\phi}g_{ik}
-\nabla_i {\phi} \nabla_l {\phi}_l g_{jk} -\nabla_j {\phi}\nabla_k
{\phi}g_{il} +|\nabla\phi|^2g_{il}g_{jk}-
|\nabla\phi|^2g_{ik}g_{lj}]
\end{split}
\end{equation}

\par We also recall the formula (for $n> 3$)\footnote[2]{For $n=3$ the Cotton tensor
 ($C_{jkl}=\nabla_kP_{lj}-\nabla_lP_{kj}$)  is
 conformally invariant and can be thought of as a substitute for the Weyl tensor.}:

\begin{equation}
\label{cotton}
\nabla_aP_{bc}-\nabla_bP_{ac}=\frac{1}{n-3}
\nabla^dW_{abcd}
\end{equation}

Now, consider our Riemannian operator  $L_{g}(f)=\sum_{h\in H} a_h
C^h_g(f)$, where each contraction $C^h_g(f)$ is in the form
(\ref{contr1}). {\it Throughout the rest of this paper we will be
assuming that $L_g(f)$ is conformally invariant of bi-degree
$(w,\kappa\cdot w-K)$}. We will re-write $L_g(f)$ as a linear
combination of contractions involving factors $\nabla^{(a)}f$,
differentiated Weyl tensors
 and also tensors of the form
 $S\nabla^{(p)}_{r_1\dots r_p}P_{ab}$ (which stands for the
fully symmetrized part of the $(p+2)$-tensor
$\nabla^{(p)}_{r_1\dots r_p}P_{ab}$). Explicitly, our ``new''
 complete contractions will be in the form\footnote[3]{For $n=3$,
we recall the discussion after theorem \ref{careful}. Thus, if
$L_g(f)$ contains complete contractions with curvature factors we
derive that $L^\sigma_g(f)=\sum_{b\in B_1\bigcup B_2} a_b
C^b_g(f)$, where the contractions indexed in $B_1$ are in the form
$contr(S\nabla^{(p)}P\otimes\nabla^{(y)}f)$ and the ones indexed
in $B_2$ must be in the form
$contr(\nabla^l\nabla^{(m)}C_{jkl}\otimes\nabla^{(p)}f)$.}

\begin{equation}
\label{Wcontr}
contr(\nabla^{(m_1)}W\otimes\dots\otimes\nabla^{(m_s)}W\otimes
S\nabla^{(p_1)}P\otimes\dots\otimes S\nabla^{(p_r)}P\otimes
\nabla^{(a_1)}f\otimes\dots\otimes\nabla^{(a_\kappa)}f)
\end{equation}

\par This can be done easily: Starting with  any complete
contraction $C^h_g(f)$ in $L_g(f)$ (in the form (\ref{contr1})), we only have to decompose the curvature tensor
as in (\ref{weyl}) and then repeatedly apply the equations
 (\ref{cotton}) and (\ref{curvature}), to express $C^h_g(f)$ as a linear combination
 of contractions in the form (\ref{Wcontr}). Thus,  we re-write $L_g(f)$:

\begin{equation}
\label{chomeur} L_{g}(f)=\sum_{u\in U} a_u C^u_{g}(f)
\end{equation}
where each $C^u_{g}(f), u\in U$ is in the form (\ref{Wcontr}).
\newline

\par In order to state our Propositions below, we need a final
piece of notation:

Consider any complete contraction $C_{g}(f)$ in the form
(\ref{Wcontr}), or even more generally in the form:

\begin{equation}
\label{Wcontr'}
contr(\nabla^{(m_1)}W\otimes\dots\otimes\nabla^{(m_s)}W\otimes
S\nabla^{(p_1)}P\otimes\dots\otimes
S\nabla^{(p_r)}P\otimes\nabla^{(\nu_1)}R\otimes\dots\otimes\nabla^{(\nu_t)}R\otimes
\nabla^{(a_1)}f\otimes\dots\otimes\nabla^{(a_\kappa)}f)
\end{equation}
For  any factor $F_h=\nabla^{(p)}f$ or $F_h=\nabla^\nu R$, we define
$\beta[F_h]$, $\gamma[F_h]$ as in the case of complete
contractions in the form (\ref{contr1}). For any factor $F_h$ in
the form $\nabla^{(m)}W_{ijkl}$, we define
 $\gamma[F_h]=\iota[F_h]-\tau[F_h]-2$. Also, for any factor $F_h$ of the
 form $S\nabla^{(p)}P$, we define $\gamma[F_h]=\iota[F_h]-\tau[F_h]-1$.
(Recall that $\iota[F_h]$ stands for the total number of indices in $F_h$
 and $\tau[F_h]$ stands for the total number of internal contractions in $F_h$).

\par In general, for any complete contraction $C_{g}(f)$ in the
form (\ref{Wcontr}) or (\ref{Wcontr'}), we define
$\beta[C_{g}(f)]$ to be the maximum among the numbers $\beta[F_h]$,
for factors $F_h$ in the form $\nabla^{(p)}f$. We also define
$\gamma[C_{g}(f)]$ to be the maximum of the numbers $\gamma[F_h]$,
where $F_h$ can be any factor in $C_g(f)$.

\par A technical tool that will be useful further down is the following:

\begin{observation}
\label{sochi} Consider a Riemannian differential operator
$L_{g}(f)$, expressed in the form:

\begin{equation}
\label{theform} L_{g}(f)=\sum_{h\in H} a_h C^h_{g}(f)
\end{equation}
where each $C^h_{g}(f)$ is in the form (\ref{contr1}). Suppose
that for each $h\in H$ we have $\beta[C^h_{g}(f)]\le \tau_1$,
$\gamma[C^h_{g}(f)]\le \tau_2$.

 Then, as explained above, we write $L_{g}(f)$ as a linear
 combination:

\begin{equation}
\label{theform2} L_{g}(f)=\sum_{h\in H'} a_h C^h_{g}(f)
\end{equation}
where each $C^h_{g}(f)$ is in the form (\ref{Wcontr}). It follows
that for each $h\in H'$ we have $\beta[C^h_{g}(f)]\le \tau_1$,
$\gamma[C^h_{g}(f)]\le \tau_2$. The converse is also true.
\end{observation}
In view of the above, we may consider a Riemannian operator and
write it as a linear combination of contractions in either of the
forms (\ref{contr1}), (\ref{Wcontr}) and unambiguously say
$\beta[L_{g}(f)]\le \tau_1$ and/or $\gamma[L_{g}(f)]\le \tau_2$.
Thus, we see that since $L_g(f)$ fulfils the extra restrictions
(when they are applicable) when written as a linear combination of
contractions in the form (\ref{contr1}), it still fulfils the
extra restrictions when written as a linear combination of
contractions in the form (\ref{Wcontr}).
\newline

Now, to prove our Theorems \ref{ayto}, \ref{ayto2}, \ref{careful}
we consider $L_{g}(f)=\sum_{u\in U} a_u C^u_g(f)$, (each
$C^u_g(f)$ in the form (\ref{Wcontr})), and we break the index set
$U$ into subsets:

\par Firstly, we pick out the complete contractions $C^u_{g}(f)$
(in the form (\ref{Wcontr})) with the minimum
 number of factors, say $\sigma$.

\par We index the contractions with $\sigma$ factors in
$U_\sigma\subset U$. We then further subdivide $U_\sigma$ into
subsets $U_{\sigma,a}, a=0,1,\dots
 ,\sigma-\kappa$ according to the number of factors
$S\nabla^{(p)}P$ in
 $C^u_{g}(f)$: We say $u\in U_{\sigma,a}$ if $C^u_{g}(f)$
(in the form (\ref{Wcontr})) contains $a$ factors
$S\nabla^{(p)}P$.

\par Our first Proposition is then the following:

\begin{proposition}
\label{queen} Consider the maximum $a$ for which $U_{\sigma,a}\ne
\emptyset$, and denote it by $a_M$. Suppose $a_M>0$. We then claim
that we can write:

\begin{equation}
\label{piece1} \sum_{u\in U_{\sigma,a_M}} a_u C^u_{g}(f)= \sum_{j\in
J} a_j C^j_{g}(f)
\end{equation}
where each $C^j_{g}(f)$ is a complete contraction in the
 form (\ref{contr1}) with at least $\sigma+1$ factors.
Furthermore each $C^j_{g}(\phi)$ satisfies the extra restrictions
$\beta[C^j_{g}(\phi)]< k$ and $\gamma[C^j_{g}(\phi)] <
\frac{n}{2}$, whenever these restrictions are applicable.
\end{proposition}

\par We observe that if we can prove the above we may just replace the sublinear combination
$\sum_{u\in U_{\sigma,a_M}} a_u C^u_{g}(f)$ in $L_{g}(f)$ by the
right hand side of (\ref{piece1}). Thus, if we can prove
Proposition \ref{queen}, by iterative repetition we
 reduce ourselves to showing our theorems under the additional
 assumption that each $C^u_{g}(f)$, $u\in U_\sigma$ has no
factors $S\nabla^{(p)}P$.

\par Under this assumption, we define $U^{0,*}_\sigma\subset
 U^0_\sigma$ to stand for the index set of complete
 contractions in the form (\ref{Wcontr}) with length $\sigma$,
no factors $S\nabla^{(p)}P$ and also no internal contractions
among any of its factors. We claim:

\begin{proposition}
\label{queen2} Suppose that the maximum $a$ for which $U_{\sigma,a}
\ne \emptyset$ is 0. Consider the sublinear combination
$\sum_{u\in U_{\sigma,0}^{*}} a_u C^u_{g}(f)$; we claim that we can
construct a Weyl operator
 $\sum_{u\in U_{\sigma,0}^{*}} a_u
\tilde{C}^u_{\tilde{g}}(\tilde{u}_w)$, so that:

\begin{equation}
\label{piece2} \sum_{u\in U^{0,*}_\sigma} a_u C^u_{g}(f)-
\sum_{u\in U^{0,*}_\sigma} a_u
\tilde{C}^u_{\tilde{g}}(\tilde{u}_w)=\sum_{u\in U'} a_u
C^u_{g}(f)+\sum_{j\in J} a_j C^j_{g}(f)
\end{equation}
Here each $C^u_{g}(f)$, $u\in U'$, is in the form (\ref{Wcontr})
with $\sigma$ factors, each in the form $\nabla^{(m)}W$ or
$\nabla^{(p)}f$ and at least one of which has an internal
contraction. The contractions $C^j_{g}(f)$ are as in the previous
Proposition.
 Moreover, the contractions on the RHS satisfy the extra
 restrictions whenever they are applicable.
\end{proposition}

\par Clearly, if we can show the above we will be reduced to
showing our Theorems in the case where each $C^u_{g}(f)$ with
$\sigma$ factors in $L_{g}(f)$ has no factors $S\nabla^{(p)}P$ and
also has at least one
 internal contraction. Under that assumption, we
then break up the index set $U_\sigma$ into subsets
$U^\delta_\sigma$, according to the rule that $u\in
U^\delta_\sigma$ if and only if $C^u_{g}(f)$ has $\delta$ internal
contractions. We then have our last and hardest claim:

\begin{proposition}
\label{queen3} Suppose that all the contractions with $\sigma$
factors in $L_{g}(f)=\sum_{u\in U} a_u C^u_{g}(f)$ are in the form
(\ref{Wcontr}) with no factors $S\nabla^{(p)}P$.
 Suppose also that the minimum $\delta$
 for which $U^\delta_\sigma\ne\emptyset$ is $\mu>0$.
 We claim that we can write:

\begin{equation}
\label{piece3} \sum_{u\in U^{\mu}_\sigma} a_u C^u_{g}(f)=
\sum_{u\in U'} a_u C^u_{g}(f)+\sum_{j\in J} a_j C^j_{g}(f)
\end{equation}
Here each $C^u_{g}(f)$, $u\in U'$, is in the form (\ref{Wcontr})
with $\sigma$ factors, each in the form $\nabla^{(m)}W$ or
$\nabla^{(p)}f$
 and with at least $\mu+1$ internal contractions.
The contractions on the RHS satisfy the extra restrictions
whenever they are applicable. The linear combination indexed in
$J$ is as in Proposition \ref{queen}.
\end{proposition}

\par Clearly, there is an obvious upper bound on the number of
internal contractions for any complete contraction of weight
 $-K$ (for example $K$). Therefore, if we can show the above
 Lemma then by iterative repetition we will be reduced to the
  case where each complete contraction in $L_{g}(f)$ has at
   least $\sigma+1$ factors.

\par Now, also observe that there is an obvious
upper bound on the total  number of factors for any complete
contraction of weight $-K$ and $f$-homogeneity $\kappa$ (say
$K+\kappa$). Therefore, if we can prove the above three
Propositions, by
 iterative repetition we derive our Theorems \ref{ayto},
 \ref{ayto2}, \ref{careful}\footnote[4]{For the case $n=3$ with $\kappa=1$, in the notation
  of the footnote 3 on page 19, the argument of Proposition \ref{queen3}
  will show that $\sum_{b\in B_2} a_b C^b_g(f)=0$
   modulo contractions with three factors. In conjunction with Proposition \ref{queen}
   (which holds as stated for $n=3$ and will show that
   $\sum_{b\in B_1} a_b C^b_g(f)=0$ modulo complete contractions with 3 factors),
    that will prove theorem \ref{careful} in
   this  case.}.

\subsection{Proof of the Propositions \ref{queen}, \ref{queen2}.}

{\it General discussion:} Our proof of Proposition \ref{queen}
will rely on a simple study of the transformation laws of the
complete contractions involved, under conformal changes of the
underlying metric.

\begin{definition}
\label{robers}  Given any  pair of numbers $(a,b)$, any Riemannian
operator $L_g(f)$ and any $\phi\in C^\infty (M)$ we define
$Im^{Z|(a,b)}_\phi[L_{g}(f)]$ as follows:

\begin{equation}
\label{ramsey} Im^{Z|(a,b)}_\phi[L_{g}(f)]=
\frac{\partial^Z}{\partial\lambda^Z}|_{\lambda=0}
\{e^{-b\lambda\phi}C_{e^{2\lambda\phi}g}(e^{a\lambda\phi }f) \}
\end{equation}
\end{definition}

We straightforwardly observe that $Im^{Z|(a,b)}_\phi[L_{g}(f)]$ is
just the linear combination of summands in
$e^{-b\lambda\phi}C_{e^{2\lambda\phi}g}(e^{a\lambda\phi }f)$ that
are homogenous of degree $Z$ in the function $\phi$.
\newline

\par We make an important observation that will be used often below:
Consider any operator $L_{g}(f)=\sum_{u\in U} a_u C^u_{g}(f)$ {\it
which is conformally invariant of bi-degree $(a,b)$}.  Then, for
any $Z\ge 1$ {\it and any function $\phi\in C^\infty (M)$, we must
have}:

\begin{equation}
\label{constantine} (Im^{Z|(a,b)}_\phi[L_g(f)]=)Im^{Z|(a,b)}_\phi[\sum_{u\in U} a_u
C^u_{g}(f)]=0
\end{equation}
This essentially just follows from the definition of conformal
invariance.
\newline

{\it Proof of Proposition \ref{queen}:}
\newline

\par We clarify what we will prove: Consider any
manifold $(M,g)$ and any $f\in C^\infty (M)$ and chose any point
$x_0\in M$. Then at $x_0$ the equation (\ref{piece1}) will hold.

\par In the notation of Proposition \ref{queen}, we set
$Z=a_M$. For each $u\in U^{a_M}$ (see the statement of Proposition
\ref{queen}) we denote by $C^u_{g}(f, \phi^{a_M})$ the complete
contraction which arises from $C^u_{g}(f)$ by replacing each of
the $a_M$ factors $S\nabla^{(p)}_{r_1\dots r_p}P_{r_{p+1}r_{p+2}}$
by $S\nabla^{p+2}_{r_1\dots r_{p+2}}\phi$.

\par Then, by virtue of the transformation law
(\ref{wschtrans}) and the conformal invariance of the Weyl tensor
we derive that:

\begin{equation}
\label{vanvleck} (0=)Im^{a_M|(w,w\cdot \kappa-K)}_\phi[\sum_{u\in
U} a_u C^u_{g}(f)]=(-1)^{a_M}\sum_{u\in U^{a_M}} a_u
C^u_{g}(f,\phi^{a_M})+\sum_{j\in J} a_j C^j_g(f,\phi)
\end{equation}
here each of the contractions $C^j_g(f,\phi)$ is in the general
form:
$$contr(\nabla^{(m_1)}R\otimes\dots\otimes
\nabla^{m_s}R\otimes\nabla^{(p_1)}f\otimes\dots \otimes
\nabla^{(p_\kappa)}f\otimes\nabla^{(y_1)}\phi\otimes\dots\otimes\nabla^{(y_{a_M})}\phi)$$
and has at  least $\sigma+1$ factors. Note that this equation
holds {\it for any $\phi\in C^\infty(M)$}.

\par Then, we just pick a function $\phi\in C^\infty(M)$ so that at our chosen $x_0\in M$ we
have that for every $p>1$:
$$S\nabla^{(p)}_{r_1\dots r_p}\phi(x_0)=-S\nabla^{p-2}_{r_1\dots
r_{p-2}}P_{r_{p-1}r_p}(g)(x_0)$$
 while if $p\le 1$ we have  $S\nabla^{(p)}\phi(x_0)=0$.

 \par For this value of the function $\phi$ it follows that
 (\ref{vanvleck}) implies Proposition \ref{queen}
when the extra restrictions are not applicable. When the
extra restrictions are applicable, we must also observe that
the correction terms (with more factors) that we are
 introducing in (\ref{vanvleck})
also satisfy the extra restrictions-this follows by the same
arguments as in the proof of Lemma \ref{appendix} (see the appendix).
$\Box$
\newline

{\it Proof of Proposition \ref{queen2}:}
\newline

\par We start by recalling that under the hypothesis of Proposition \ref{queen2},
 all the complete contractions $C^u_{g}(f)$
with $\sigma$ factors that appear in the expression for $L_{g}(f)$
must have $\sigma-\kappa$ factors $\nabla^{(m)}W$. In other words,
they will be in the form:

\begin{equation}
\label{onlyW} contr(\nabla^{(m_1)}_{r_1\dots
r_{m_1}}W_{ijkl}\otimes\dots \otimes
\nabla^{(m_{\sigma-\kappa})}_{u_1\dots u_{m_{\sigma-\kappa}}}
W_{i'j'k'l'}\otimes \nabla^{(p_1)}_{y_1\dots y_{p_1}}f\otimes
\dots \otimes\nabla^{(p_\kappa)}_{t_1\dots t_{p_\kappa}}f)
\end{equation}

\par For each complete contraction $C^l_{g}(f)$ in the above form with no internal contractions
we will construct a complete contraction in the ambient metric,
$C^l_{{\tilde{g}}^{n+2}} (\tilde{u}_w)$:
\begin{equation}
\label{ambient} contr({\tilde{\nabla}}_{r_1\dots
r_{m_1}}^{(m_1)}{\tilde{R}}_{i_1j_1k_1l_1} \otimes \dots \otimes
{\tilde{\nabla}}_{v_1\dots
v_{m_{\sigma-\kappa}}}^{(m_{\sigma-\kappa})}{\tilde{R}}_{i_sj_sk_sl_s}\otimes\tilde
{\nabla}^{(p_1)}_{y_1\dots
y_{p_1}}\tilde{u}_w\otimes\dots\otimes\tilde
{\nabla}^{(p_\kappa)}_{t_1\dots t_{p_\kappa}}\tilde{u}_w)
\end{equation}
 which is obtained from $C^l_{g}(f)$ by just replacing each factor $\nabla^{(m)}W$ by
  a factor $\tilde{\nabla}^{(m)}\tilde{R}$
 and each factor $\nabla^{(p)}f$ by $\tilde{\nabla}^{(p)}\tilde{u}_w$
and then performing the same contractions (only with respect
 to the metric $\tilde{g}$).
 We now show that this Weyl operator
$\sum_{u\in U^{0,*}_\sigma} a_u
\tilde{C}^u_{\tilde{g}}(\tilde{u})$ satisfies the conclusion of
Proposition \ref{queen2}.
\newline

{\it Proof:} In order to see this claim, let us recall some
notation from subsection \ref{amb}. We start with $(M,g)$ and
$f\in C^\infty (M)$ and we perform the ambient metric construction
picking some $x_0\in M$ and mapping it to $\tilde{x}_0=(1,x_0,0)$
in $\tilde{G}$. Recall that if the coordinates of $(M,g)$ are
$\{x^1,\dots ,x^n\}$, then there is a special coordinate system
for the ambient manifold $(\tilde{G}^{n+2},\tilde{g})$ of
 the form:
$\{t=x^0,x^1,\dots ,x^n,\rho=x^{n+1}\}$.

\par Now, let us furthermore recall the form of the ambient metric on
$G^{n+1}\subset \tilde{G}$. In the coordinate system $\{x^0,\dots
, x^{n+1}\}$ the ambient metric at $\tilde{x}_0$ is of the form:

\begin{equation}
\label{form}\tilde{g}^{n+2}_{ab}dx^adx^b =2dx^0dx^{n+1}+
\sum_{i,j=1}^ng_{ij}dx^idx^j
\end{equation}
where $0\le a,b\le n+1$ and $1\le i,j\le n$. We denote by
$X^0,X^1,\dots ,X^n,X^\infty$ the vector fields that correspond to
the coordinates $\{ x^0,x^1,\dots ,x^n,x^{n+1}\}$. In view of the
form (\ref{form}) of the ambient metric in this coordinate system,
we observe that for each pair of indices ${}_a,{}_b$ in any
$\tilde{C}^u_{\tilde{g}}(\tilde{u}_w)$ that are contracting
against each other, if we assign the value ${}_\infty$ or ${}_0$
to one of the indices, then we must assign the value ${}_0$ or
${}_\infty$ to the other.

\par Now, in order to express a complete contraction in the form
(\ref{ambient}) as a linear combination of complete contractions
in the form (\ref{Wcontr}), we will have to express the components
of each tensor $\tilde{\nabla}_{r_1\dots
r_m}^{(m)}\tilde{R}_{ijkl}(\tilde{g})$ and each tensor
$\tilde{\nabla}^{(p)}_{r_1\dots r_p}\tilde{u}_w(\tilde{g})$ in
terms of the tensors $\nabla^{(m)}W(g)$, $\nabla^{(p)}P(g)$ and
$\nabla^{(y)}f(g)$.

\par Using the Christoffel symbols of $\tilde{g}$ with respect to
the special coordinate system we can see the following: Consider
any component $T_{r_1\dots r_{m+4}}=\tilde{\nabla}^{(m)}_{r_1\dots
r_m}\tilde{R}_{r_{m+1}\dots r_{m+4}}(\tilde{g})$ with $\delta$ of
the indices ${}_{r_1},\dots ,{}_{r_{m+4}}$ being $\infty$'s,
$\epsilon$ being $0$'s and the rest having values between $1$ and
$n$. Let us suppose that the indices that have values between $1$
and $n$ are precisely $r_{a_1},\dots r_{a_q}$. It then follows
from standard computations
 on the ambient metric that:

\begin{equation}
\label{arkansas} T_{r_{1}\dots
r_{m+4}}=\sum_{h=0}^{\frac{m+4-\delta}{2}}F^h_{r_{a_1}\dots
r_{a_q}}+\sum_{j\in J} a_j F^j_{r_{a_1}\dots r_{a_q}}
\end{equation}
where each $F^h_{r_1\dots r_{m+4}}$ stands for a linear
combination of tensor products in the form $\nabla^{(m+h)}W\otimes
g\otimes\dots\otimes g$, where the factor $\nabla^{(m+h)}W$ has
$\delta+h$ internal contractions (thus, we observe that
$\gamma[F^h]\le \gamma[C^u_{g}(f)]$). If $\delta=\epsilon=0$ then
$F^0_{r_1\dots r_{m+4}}=\nabla^{(m)}_{r_1\dots r_m}W_{r_{m+1}\dots
r_{m+4}}$.

 The tensors $F^j_{r_1\dots r_{m+4}}$
stand for linear combinations of tensor products of the form
$\nabla^{(m_1)}R\otimes\dots\otimes\nabla^{(m_u)}R\otimes
g\otimes\dots\otimes g$ ($u\ge 2$), where each factor
$\nabla^{(m_y)}R$ satisfies $\gamma[\nabla^{(m_y)}R]\le
\gamma[C^u_{g}(f)]$.

\par By complete analogy, consider any component $T_{r_1\dots
r_{p}}=\tilde{\nabla}^{(p)}_{r_1\dots r_p}\tilde{u}_w(\tilde{g})$
with $\delta$ of the indices ${}_{r_1},\dots ,{}_{r_{m+4}}$ being
$\infty$'s, $\epsilon$ being $0$'s and the rest having values
between $1$ and $n$. Let us suppose that the indices that have
values between $1$ and $n$ are precisely ${}_{r_{a_1}},\dots ,
{}_{r_{a_q}}$. It then follows (from standard computations
 on the ambient metric) that:

\begin{equation}
\label{arkansas2} T_{r_1\dots
r_{p}}=\sum_{h=0}^{\frac{p-\delta}{2}}F^h_{r_{a_1}\dots
r_{a_q}}+\sum_{j\in J} a_j F^j_{r_{a_1}\dots r_{a_q}}
\end{equation}
where each $F^h_{r_1\dots r_p}$ stands for a linear combination of
tensor products in the form $\nabla^{(p+h)}f\otimes
g\otimes\dots\otimes g$, where the factor $\nabla^{(p+h)}f$ has
$\delta+h$ internal contractions (thus, we observe that
$\gamma[F^h]\le\gamma[C^u_{g}(f)]$ and
$\beta[F^h]\le\beta[C^u_{g}(f)]$). If $\delta=\epsilon=0$ then
$F^0_{r_1\dots r_p}=\nabla^{(p)}_{r_1\dots r_p}f$.

 The tensors $F^j_{r_1\dots r_{m+4}}$
stand for linear combinations of tensor products of the form
$\nabla^{(p')}f\otimes
\nabla^{(m_1)}R\otimes\dots\otimes\nabla^{(m_u)}R\otimes
g\otimes\dots\otimes g$ ($u\ge 1$), where we have
$\gamma[\nabla^{(p')}f]<\gamma[C^u_{g}(f)]$,
$\beta[\nabla^{(p')}f]<\beta[C^u_{g}(f)]$,
$\gamma[\nabla^{(m_y)}R]<\gamma[C^u_{g}(f)]$. Therefore, replacing
the factors of $\tilde{C}^u_{\tilde{g}}(\tilde{u}_w)$ by the right
hand sides of (\ref{arkansas}), (\ref{arkansas2}) we derive our
claim. $\Box$

\section{Proof of Proposition \ref{queen3}.}

{\it Outline:} This (somewhat technical) section is divided in
three parts: In \ref{bidfor} we recall some facts about identities
holding ``formally'' which will be needed in the proof. In
\ref{proof} we introduce some notation to claim Proposition
\ref{finn}: Very roughly, we consider the sublinear combination
$\sum_{u\in U^\mu_\sigma} a_u C^u_g(f)$ in the statement of
Proposition \ref{queen3}, and for each $u\in U^\mu_\sigma$ we
formally construct new complete contractions
$C_g^{u,\iota|i_1\dots
i_\mu}(f)\nabla_{i_1}\upsilon\dots\nabla_{i_\mu}\upsilon$ by
formally replacing each factor $\nabla^{(m)}_{r_1\dots
r_m}W_{ijkl}$ by $\nabla^{(m)}_{r_1\dots r_m}R_{ijkl}$ (times a
constant) and also replacing each {\it internal contraction}
$(\nabla^a,{}_a)$ by an expression $(\nabla^a\upsilon,{}_a)$.
Proposition \ref{finn} then claims that at the linearized level
$\sum_{u\in U^\mu_\sigma} a_u C^{u,\iota|i_1\dots
i_\mu}\nabla_{i_1}\upsilon\dots\nabla_{i_\mu}\upsilon=0$. We then
show that Proposition \ref{finn} implies Proposition \ref{queen3}
(via another formal operation $Weylify[\dots]$ which formally
reverses the previous construction).

\par In subsection \ref{provefinn} we prove Proposition
\ref{finn}: Briefly, the idea is to consider the sublinear
combination  in the equation
$e^{(w\cdot\kappa-K)\phi}L_{{e^2\phi}g}(e^{w\cdot \phi})-L_g(f)=0$
that is linear in $\phi$, thus obtaining a new equation,
$Im^{1|(w,w\cdot\kappa-K)}_\phi[L_g(f)]=0$. We then derive (after
a careful study of transformation laws under conformal re-scaling)
that the sublinear combination
$Im^{1||(w,w\cdot\kappa-K)|*}_\phi[L_g(f)]$ of terms with
$\sigma+1$ factors, $\mu-1$ internal contractions and a factor
$\nabla\phi$ arises {\it only} from $\sum_{u\in U^\mu_\sigma} a_u
C^u_g(f)$ by
 replacing an internal contraction $(\nabla^a,{}_a)$ by
$(\nabla^a\phi,{}_a)$ times a constant. We then observe that
$Im^{1|(w,w\cdot\kappa-K)|*}_\phi[L_g(f)]=0$, thus deriving
equation (\ref{lyapunov5}). Lemma \ref{simp2} then essentially
completes the proof of Proposition \ref{finn}, modulo checking
that the constant referred to above is non-zero.

\subsection{Identities holding formally vs. by substitution.}
\label{bidfor}

\par We here very briefly explain the theorem B.3 in
\cite{beg:itccg} and its straightforward generalization that
appears in \cite{a:dgciI}. Theorem B.3 in \cite{beg:itccg} is
itself an extension of the work of Weyl in \cite{w:cg}.

\par This theorem deals with complete contractions involving
tensors with certain symmetries and anti-symmetries. In the case
at hand, we form complete contractions of tensor products
involving  {\it symmetric tensors} $\{T^\alpha \}_{\alpha\in A}$
(that belong to a family $A$) and {\it linearized curvature
tensors}, $\{R\}=\{R_{ijkl},R_{ijkl,r_1},\dots, R_{ijkl,r_1\dots
r_m},\dots\}$. The latter model (at the linearized level) the
symmetries and anti-symmetries of the curvature and its covariant
derivatives (see \cite{a:dgciI} for more details).

\par We can form complete contractions in
the above objects:
\begin{equation}
\label{fai} C(T,R)=contr(R_{ijkl,r_1\dots r_m}\otimes\dots
R_{i'j'k'l',t_1\dots t_u}\otimes T^{\alpha_1}_{y_1\dots
y_b}\otimes T^{\alpha_\tau}_{z_1\dots z_c})
\end{equation}
 (say with $\rho$ factors in the form
 $R_{ijkl,r_1\dots r_m}$ and $\tau$ factors
$T^\alpha_{y_1\dots y_b}$, $b\ge 1$, $\alpha\in A$; the
contractions are taken with respect to $\delta_{ij}$)
 and consider linear combinations thereof.

 There are then two notions of an identity  holding between linear
 combinations of such complete contractions: Following
\cite{beg:itccg}, we say an identity
 holds {\it by substitution} if it holds for all possible assignments of
 values to the tensors in $\{T^\alpha\}_{\alpha\in A}$ and
$\{R$\}, which satisfy the symmetry and
 anti-symmetry restrictions we have imposed. We say an identity
 holds {\it formally} if we can just prove it by virtue of applying the
 symmetries and anti-symmetries and also the distributive rule (see
 \cite{a:dgciI} for  a precise discussion). Clearly, if an identity
 holds  formally, it will then also hold by substitution. Theorem 2
 in \cite{a:dgciI} says that the converse is also true, subject to
one restriction:
\begin{proposition}
\label{likeit}
 Suppose that $C^l(T,R)$, $l\in L$ are complete contractions in the form (\ref{fai}), and the identity:

\begin{equation}
\label{kahramanoglou} \sum_{l\in L} a_l C^l(T,R)=0
\end{equation}
  holds by substitution. Moreover suppose that each $C^l(T,R)$ satisfies
$\tau+2\rho\le n$ (see the notation after (\ref{fai})). It then
follows that (\ref{kahramanoglou}) also holds formally.
\newline

\par Alternatively, suppose that for each $C^l(T,R)$ above we have
$\tau+2\rho=n+1$, but also that each $C^l(T,R)$ has one factor
$T^{\alpha_a}_x$ with rank 1, this factor only appearing once in
each contraction. Then (\ref{kahramanoglou}) again holds formally
\end{proposition}

\par Finally, let us recall the notion of {\it linearization} for
complete contractions in the form (\ref{contr1}): For any such
complete contraction $C^l_{g}(f)$ we define $lin\{C^l_{g}(f) \}$
to stand for the complete contraction in the form (\ref{fai}) that
is constructed out of $C^l_{g}(f)$ by substituting each
``genuine'' curvature term $\nabla^{(m)}_{r_1\dots r_m}R_{ijkl}$
by a term $R_{ijkl,r_1\dots r_m}$ and each ``genuine'' covariant
derivative of $f$, $\nabla^{(p)}_{r_1\dots r_p}f$ by a symmetric
 tensor $T_{r_1\dots r_p}$. We recall the following fact from
\cite{a:dgciI}:

\begin{proposition}
\label{equiv} Consider a linear combination $\sum_{l\in L} a_l
C^l_{g}(f)$ of complete contractions in the form (\ref{contr1}),
each with $\sigma$ factors. Suppose that for any $g, f$ we have an
 identity:

\begin{equation}
\label{lohkamp} \sum_{l\in L} a_l C^l_{g}(f)=\sum_{h\in H} a_h
C^h_{g}(f)
\end{equation}
 where the RHS stands for some linear
combination of complete contractions in the form (\ref{contr1})
with at least $\sigma +1$ factors. Suppose that for each $l\in L$
 we have that $lin\{ C^l_{g}(f)\}$ satisfies $\tau+2\rho\le n$.

Then we also have an identity: \begin{equation} \label{lohkamp2}
\sum_{l\in L} a_l lin\{C^l_{g}(f)\}=0
\end{equation}
 which holds formally.

 \par Alternatively, if we assume (\ref{lohkamp}) where
 each of the contractions $C^l_g(f)$ is in
 the form (\ref{contr1}) but also has an additional factor $\nabla\phi$
 (thus $C^l_g(f)$ has $\sigma+1$ factors in total)
 and also each $C^h_g(f)$ also has a factor
 $\nabla^{(y)}\phi$ (thus $C^h_g(f)$ has $\sigma+2$ factors in total)
  and furthermore $\tau+2\rho\le n+1$; then
  the linearized equation (\ref{lohkamp2}) will hold formally.
\end{proposition}

\subsection{Decomposing $\nabla^{(m)}W$ and reducing Proposition \ref{queen3} to Proposition \ref{finn}.}
\label{proof}

{\it Convention:}  For each tensor $T=\nabla^{(m)}_{r_1\dots
r_m}W_{ijkl}$, we will call the indices $i,j,k,l$ the {\it
internal indices} of $T$.
\newline

\par We will now put down some well-known identities that will be useful for
our discussion. We will be considering complete
 contractions that involve the Weyl tensor and write each such
complete contraction as a linear combination of contractions
 involving only the curvature tensor. We do this via the
formula (\ref{weyl}). It will be useful further down to be more
precise about this decomposition, when we consider
$\nabla^{(m)}W_{ijkl}$ (i.e. an iterated covariant derivative of
the Weyl tensor).

 Consider any tensor
$T=\nabla^{r_{a_1}\dots r_{a_x}} \nabla^{(m)}_{r_1\dots r_m}W_{ijkl}$
where each index ${}^{r_{a_s}}$ is contracting against the
(derivative) index ${}_{r_{a_s}}$, and all the other indices are
free (so we have an $(m+4-x)$-tensor). By (\ref{weyl}) it follows that:

\begin{equation}
\label{decompo1} \nabla^{r_{a_1}\dots r_{a_x}}
\nabla^{(m)}_{r_1\dots r_m}W_{ijkl}=\nabla^{r_{a_1}\dots
r_{a_x}}\nabla^{(m)}_{r_1\dots r_m}R_{ijkl}+\sum_{z\in
Z^{\delta=x+1}} a_z T^z(g)+ \sum_{z\in Z^{\delta=x+2}} a_z T^z(g)
\end{equation}
where $\sum_{z\in Z^{\delta= x+1}} a_z T^x(g)$ stands for a
linear combination of tensor products of the form
$\nabla^{r_{a_1}\dots r_{a_x}}\nabla^{(m)}_{r_1\dots
r_m}Ric_{sq}\otimes g_{vb}$ in the same free indices as $T$, with
the feature that there are a total of $x+1$ internal contractions
in the tensor $\nabla^{(m)}Ric_{sq}$ (including the one in the tensor
$Ric_{sq}={R^a}_{saq}$ itself). $\sum_{z\in Z^{\delta= x+2}} a_z
T^x(g)$ stands for a linear combination of tensor products of the
form $\nabla^{r_{a_1}\dots r_{a_x}}\nabla^{(m)}_{r_1\dots r_m}S\otimes
g_{vb}\otimes g_{hj}$ ($S$ is the scalar curvature)
 in the
same free indices as $T$, with the feature that there are a total
of $x+2$ internal contractions in the tensor $\nabla^{(m)}S$
(including the two in the factor $S={R^{ab}}_{ab}$ itself).

\par Now, we consider a factor $T$ in the form
$T=\nabla^{r_{a_1}\dots r_{a_x}}\nabla^{(m)}_{r_1\dots
r_m}W_{r_{m+1}r_{m+2}r_{m+3}r_{m+4}}$ where again each upper
 index ${}^{r_{a_v}}$ is contracting against the lower index
${}_{r_{a_v}}$, and moreover at least one of the indices
${}^{r_{a_v}}$ is contracting against one of the internal
 indices ${}_{r_{m+1}},\dots ,{}_{r_{m+4}}$. By applying
 (\ref{weyl}) it can be seen that:

\begin{equation}
\label{decompo2} \begin{split} &T=\nabla^{r_{a_1}\dots
r_{a_x}}\nabla^{(m)}_{r_1\dots
r_m}W_{r_{m+1}r_{m+2}r_{m+3}r_{m+4}}=
\frac{n-3}{n-2}\nabla^{r_{a_1}\dots r_{a_x}}\nabla^{(m)}_{r_1\dots
r_m} R_{r_{m+1}r_{m+2}r_{m+3}r_{m+4}} \\&+\sum_{z\in Z^{\delta=x}}
a_z T^z(g)+\sum_{z\in Z^{\delta=x+1}} a_z T^z(g)
\end{split}
\end{equation} where $\sum_{z\in Z^{\delta=x+1}} a_z T^x(g)$
stands for the same linear combination as before. Now
$\sum_{z\in Z^{\delta=x}} a_z T^x(g)$ only appears in the case
where there are two indices ${}^{r_{a_b}}, {}^{r_{a_c}}$
contracting against two internal indices ${}_{r_a},{}_{r_b}$ in
$W_{ijkl}$ (and moreover the indices ${}_{r_{a_b}},{}_{r_{a_c}}$
do not belong to the same block $[ij],[kl]$). It stands for a
linear combination of tensors $\nabla^{r_{a_1}\dots
r_{a_x}}\nabla^{(m)}Ric_{ab}$ with $x$ internal contractions, and with
the extra feature that one of the indices ${}^{r_{a_1}},\dots
,{}^{r_{a_x}}$ is contracting
 against one of the indices ${}_a,{}_b$ in $Ric_{ab}$. In fact, to
 facilitate our
discussion further down, we repeatedly apply the Ricci identity
and the contracted second Bianchi identity
$2\nabla^aRic_{ab}=\nabla_bS$ to re-write the linear combination
$\sum_{z\in Z^{\delta=x}} a_z T^z(g)$ above in the form:

\begin{equation}
\label{decompo2,5} \sum_{z\in Z^{\delta=x}} a_z
T^z(g)=\sum_{z\in Z'^{\delta=x}} a_z T^z(g) +\sum_{q\in Q} a_q
T^q(g)
\end{equation}
where each $T^z(g)$, $z\in Z'^{\delta=x}$ is a factor of the
 form $\nabla^{m+x}S$ ($S$ is the scalar curvature),
with a total of $x$ internal contractions (including the two
 internal contractions in $S={R^{ij}}_{ij}$). Also, $T^q$ are
  quadratic correction terms, partial contractions of the
form $pcontr(\nabla^{(b)}R\otimes\nabla^{(c)}R)$, moreover we have
$\gamma[\nabla^{(b)}R]\le \gamma[T]$ and $\beta[\nabla^{(b)}R]\le
 \beta[T]$ and similarly for $\nabla^{(c)}R$
($T$ is the left hand side of (\ref{decompo2})).
\newline

\par Now, we will prove our proposition \ref{queen3} by virtue of the next
 proposition, for which we will need a little more notation.

\par For each contraction $C^u_{g}(f)$ in the form (\ref{onlyW}), we define
$C^{u,\iota}_{g}(f)$ to stand for the complete contraction
 (times a constant), in the form:

\begin{equation}
\label{greatman} contr(\nabla^{y_1\dots
y_b}\nabla^{(m_1)}_{r_1\dots r_m}R_{ijkl}\otimes\dots\otimes
\nabla^{w_1\dots w_z}\nabla^{(m_x)}_{t_1\dots t_{m_x}}
R_{i'j'k'l'}\otimes\nabla^{(p_1)}_{v_1\dots v_{p_1}}f\otimes
\dots\otimes \nabla^{(p_b)}_{v_1\dots v_{p_q}}f)
\end{equation}
which arises from $C^u_{g}(f)$ (in the form (\ref{onlyW})) by
replacing each factor $\nabla^{u_1\dots u_x}\nabla^{(m)}_{r_1\dots
r_m}W_{ijkl}$ with no internal
 contractions involving internal indices by a factor
$\nabla^{u_1\dots u_x}\nabla^{(m)}_{r_1\dots r_m}R_{ijkl}$ and each
factor $\nabla^{u_1\dots u_x}\nabla^{(m)}_{r_1\dots r_m}W_{ijkl}$ with
at least one internal contraction involving an internal index
by a factor $\frac{n-3}{n-2} \nabla^{u_1\dots
u_x}\nabla^{(m)}_{r_1\dots r_m}R_{ijkl}$.

\begin{definition}
\label{tereza}
 For each $u\in U^\mu_\sigma$ we consider $C^{u,\iota}_{g}(f)$
  in the form (\ref{greatman}) and we construct a pair of
lists $(L_1,L_2)$: $L_1$ stands for the list $(\delta_1,\dots
,\delta_x)$, where $\delta_i$ is the number of internal
contractions in the factor $T_i(=\nabla^{(m_i)}R)$ in
$C^{u,\iota}_{g}(f)$. $L_2$ stands for the list
$(\delta_{x+1},\dots,\delta_{x+b})$ where $\delta_i$ is against
the number of internal contractions in the factor
$T_i(=\nabla^{(p)}f)$ in $C^{u,\iota}_{g}(f)$. We then define
$RL_1,RL_2$ to stand for the decreasing rearrangements of the
lists $L_1,L_2$ after we erase the $0$-entries. We call
$(RL_1,RL_2)$ the character of $C^u_{g}(f)$ and denote it by
$\vec{\lambda}(u)$. We will also denote by $\Lambda$ the set of
all pairs of lists $\Lambda=\{\vec{\lambda}(u)\}_{u\in
U^\mu_\sigma}$.
\end{definition}

\par Accordingly, we subdivide $U^\mu_\sigma$ into subsets
$U^{\mu,\vec{\alpha}}_\sigma$, $\vec{\alpha}\in \Lambda$ according
to the rule that $u\in U^{\mu,\vec{\alpha}}_\sigma$ if and only if
$C^u_{g}(f)$ has $\vec{\lambda}(u)=\vec{\alpha}$.

\par Then, for every $u\in U^\mu_\sigma$, we define $C^{u,\iota|i_1\dots i_\mu}_{g}(f)$ to
 stand for the tensor field that arises from
$C^{u,\iota}_{g}(f)$ by replacing
 each internal contraction $(\nabla^a,{}_a)$ by a free index ${}_a$
 (in other words we erase $\nabla^a$ and make the index ${}_a$ free).
  We then construct the complete contractions:

$$C^{u,\iota|i_1\dots i_\mu}_{g}(f)\nabla_{i_1}\upsilon\dots
\nabla_{i_\mu}\upsilon$$ which arise from each
$C^{u,\iota|i_1\dots i_\mu}_{g}(f)$ by just contracting each of
the free indices against a factor
 $\nabla\upsilon$ ($\upsilon$ is some arbitrary scalar
 function). The notion of character naturally extends to tensor
 fields, or to complete contractions in
 the above form, where instead
  of internal contractions we count free
 indices or the numbers of indices that contract
  against factors $\nabla\upsilon$, respectively.

\par We then consider the linearizations of the
  complete contractions above, which we denote by:

\begin{equation}
\label{papa} lin\{ C^{u,\iota|i_1\dots i_\mu}_{g}(f)\nabla_{i_1}
\upsilon\dots\nabla_{i_\mu}\upsilon\}
\end{equation}
(we will denote by $R_{ijkl,r_1\dots r_m}$ the linearized
curvature factor that replaces $\nabla^{(m)}_{r_1\dots
r_m}R_{ijkl}$, by $\Phi_{r_1\dots r_p}$ the symmetric tensor that
replaces $\nabla^{(p)}_{r_1\dots r_p}f$, and by $\upsilon_s$ the
vector that replaces $\nabla_s\upsilon$).

\par We claim a new Proposition, which will imply
our  Proposition \ref{queen3}:

\begin{proposition}
\label{finn} In the above equation, we claim that for each
$\vec{\alpha}\in \Lambda$:

\begin{equation}
\label{korevaar} \sum_{u\in U^{\mu,\vec{\alpha}}_{\sigma}}a_u lin\{
C^{u,\iota|i_1\dots i_\mu}_{g}(f)\nabla_{i_1}
\upsilon\dots\nabla_{i_\mu}\upsilon\}=0
\end{equation}
and moreover the above holds formally.
\end{proposition}

{\it Proof that Proposition \ref{finn} implies Proposition
\ref{queen3}:}
\newline

\par Clearly, in order to prove Proposition \ref{queen3} it would
 suffice to prove that for every $\vec{\alpha}\in \Lambda$:

\begin{equation}
\label{piece1*} \sum_{u\in U^{\mu,\vec{\alpha}}_\sigma} a_u
C^u_{g}(f)= \sum_{u\in U'} a_u C^u_{g}(f)+\sum_{j\in J} a_j
C^j_{g}(f)
\end{equation}
where the right hand side is as in the statement of
 Proposition \ref{queen3}. We will show this below.
 We start by putting down a few identities.

We recall that the Weyl tensor $W_{ijkl}$ is antisymmetric in
 the indices $i,j$ and also $W_{ijkl}=W_{klij}$. It also
satisfies the first Bianchi identity. Nevertheless, it does not
satisfy the second Bianchi identity. We now present certain
substitutes for the second Bianchi identity:

Firstly, if the indices $r,i,j,k,l$ are all free we have that:

\begin{equation} \label{in&out}
\nabla_rW_{ijkl}+\nabla_jW_{rikl}+\nabla_iW_{jrkl}=\sum
(\nabla^sW_{srty}\otimes g)
\end{equation}
where the symbol $\sum (\nabla^sW_{srty}\otimes g)$ stands for a
linear combination of tensor products of the three-tensor
$\nabla^sW_{srty}$ (i.e., essentially the Cotton tensor) with an
un-contracted metric tensor. The exact form of $\sum
(\nabla^sW_{srty}\otimes g)$ is not important so we do not
 write it down.

\par On the other hand, if the indices $i,j,k,l$ are free we
 then have:

\begin{equation}
\label{in&out2}
\nabla^s_sW_{ijkl}+\frac{n-2}{n-3}\nabla^s_jW_{sikl}+\frac{n-2}{n-3}\nabla^s_iW_{jskl}=
\sum (\nabla^{st}W_{svtr}\otimes g)+Q(R)
\end{equation}
where the symbol $\sum (\nabla^{st}W_{svtr}\otimes g)$ stands
for a linear combination of tensor products:
$\nabla^{st}W_{svtr}\otimes g_{ab}$ ($g_{ab}$ is an un-contracted
metric tensor-note that there are {\it two} internal contractions
in the factor $\nabla^{ik}W_{ijkl}$). $Q(R)$ stands for a
quadratic
 expression in the curvature tensor (without covariant
 derivatives). Again the exact form of these expressions is
not important so we do not write them down.

\par Whereas, if the indices $r,i,j,l$ are free we will then
have that:

\begin{equation}
\label{in&out2c}
\nabla^k\nabla_rW_{ijkl}+\nabla^k\nabla_jW_{rikl}+\nabla^k\nabla_iW_{jrkl}=
 Q(R)
\end{equation}

\par Finally:

\begin{equation}
\label{in&out2d}
\nabla^{ri}\nabla_rW_{ijkl}+\nabla^{ri}\nabla_jW_{rikl}+\nabla^{ri}\nabla_iW_{jrkl}= Q(R)
\end{equation}

\par Of course, if we take covariant derivatives of these
equations, they continue to hold.  We will collectively call these
identities the ``fake'' second Bianchi identities.

\par A lemma that will be useful in the more technical parts of this proof is the following:

\begin{lemma}
\label{appendix} Consider any complete contraction $C_g(f)$ in
either of the forms (\ref{contr1}), (\ref{Wcontr}).
 Suppose we apply either the identity (\ref{curvature}),
or any of the fake second Bianchi identities, and thus write:
$C_g(f)=\sum_{k\in K} a_k C^k_g(f)$.

\par We then claim that if $C_g(f)$ satisfies the extra restrictions then so does each $C^k_g(f)$.
\end{lemma}

We will prove this lemma in the appendix. Let us now make a few
easy observations.

\par  Any complete contraction $C^u_{g}(f)$ in the form
(\ref{onlyW}) that has two
 antisymmetric indices ${}_i,{}_j$ or ${}_k,{}_l$ in a
given factor $T=\nabla^{(m)}_{r_1\dots r_m}W_{ijkl}$ contracting
 against two derivatives in the same factor can be written: $C^u_{g}(f)=
 \sum_{r\in R} a_r C^r_{g}(f)$
each $C^r_{g}(f)$ with $\sigma+1$ factors. This is straightforward
from the antisymmetry of the indices $i,j$ and $k,l$ and the Ricci
identity. Moreover we have that each $C^r_{g}(f)$ will satisfy the
extra restrictions when they
 are applicable (by Lemma \ref{appendix}).
  Thus we may prove our Proposition
under the extra assumption that all complete contractions
$C^u_{g}(f)$ have no factor $T=\nabla^{(m)}_{r_1\dots r_m}
W_{ijkl}$ with two antisymmetric indices ${}_i,{}_j$ or
${}_k,{}_l$ contracting against two derivative indices in $T$.
Moreover, since the Weyl tensor is trace-free we may assume that
each $C^u_{g}(f)$, $u\in U_\sigma$ has no factor $\nabla^{(m)}W$
with two internal indices contracting between themselves.
\newline

\par Now  we will show that
Proposition \ref{finn} implies Proposition \ref{queen3}:

\begin{definition}
\label{weylify}  We define an operation $Weylify$ that acts on
linearized complete contractions in the form (\ref{papa}) as
follows: We identify
 the indices in each linearized factor
$R_{ijkl,r_1\dots r_m}$ and $\Phi_{w_1\dots w_a}$ that are
contracting against factors $\upsilon_s$. Then, we pick out
each factor $\Phi_{w_1\dots w_a}$ where
 the indices $w_{h_1},\dots ,w_{h_b}$ are contracting against
  factors $\upsilon_s$ and replace it by a factor
$\nabla^{w_{h_1}\dots w_{h_b}}\nabla^a_{w_1\dots w_a}f$ (and we
also erase the factors $\upsilon_s$), thus obtaining a factor with
$h_b$ internal contractions.

\par Moreover, we pick out each factor
$R_{r_{m+1}r_{m+2}r_{m+3}r_{m+4},r_1\dots r_m}$ and we identify
the indices $r_{h_1}, \dots ,r_{h_b}$  that are contracting
against factors $\upsilon_s$. We then inquire whether any of the
internal indices $r_{m+1},\dots ,r_{m+4}$ are contracting against
a factor $\upsilon_s$. If not, and we replace the factor
$R_{ijkl,r_1\dots r_m}$ by a factor $\nabla^{r_{h_1}\dots
r_{h_b}}\nabla^{(m)}_{r_1\dots r_m} W_{ijkl}$. If there are internal
indices contracting against factors $\upsilon_s$ we then replace
 $R_{r_{m+1}r_{m+2}r_{m+3}r_{m+4},r_1\dots r_m}$ by
a factor $\frac{n-2}{n-3}\nabla^{r_{h_1}\dots r_{h_b}}
\nabla^{(m)}_{r_1\dots r_m}W_{ijkl}$. This operation $Weylify[\dots ]$
extends to linear
 combinations.
\end{definition}

\par By definition, we observe that:

\begin{equation}
\label{sinister} Weylify \{ \sum_{u\in
U^{\mu,\vec{\alpha}}_{\sigma}}a_u lin\{ C^{u,\iota|i_1\dots
i_\mu}_{g}(f)\nabla_{i_1} \upsilon\dots\nabla_{i_\mu}\upsilon\}\}=
\sum_{u\in U^{\mu,\vec{\alpha}}_{\sigma}}a_u C^{u}_{g}(f)
\end{equation}

\par  By virtue of the above and also of equations
(\ref{in&out}), (\ref{in&out2}), (\ref{in&out2c}),
(\ref{in&out2d}),  we observe that if we {\it repeat} the sequence
of permutations by which we make
 the left hand side of (\ref{korevaar}) formally zero, we
can make $\sum_{u\in U^{\mu,\vec{\alpha}}_{\sigma}}a_u
C^{u}_{g}(f)$ formally zero, modulo introducing correction terms
by virtue of the right hand sides of the equations (\ref{in&out}),
(\ref{in&out2}), (\ref{in&out2c}), (\ref{in&out2d}) and also by
virtue of the
 Ricci identity when we interchange adjacent derivative indices.

\par But by inspection of the formulas (\ref{in&out}),
(\ref{in&out2}), (\ref{in&out2c}), (\ref{in&out2d}) and also by
using Lemma \ref{appendix}, it follows that the correction terms
that arise thus are in the form:

$$\sum_{u\in U'} a_u C^u_{g}(f)+\sum_{j\in J} a_j
C^j_{g}(f)$$ (in the notation of Proposition \ref{queen3})).
$\Box$

\subsection{Proof of Proposition \ref{finn}.}
\label{provefinn}

\par Our point of departure in this proof is the equation:

$$Im^{1|(w,w\cdot\kappa -K)}_\phi[L_{g}(f)]=0.$$

\par In order to apply our arguments, we will be writing out
$Im^{1|(w,w\cdot\kappa -K)}_\phi[L_{g}(f)]$
 as a linear combination of complete contractions in
the form:

\begin{equation}
\label{skeleton} contr(\nabla^{(m_1)}R\otimes\dots\otimes
\nabla^{(m_t)}R\otimes\nabla^{(p)}\phi)
\end{equation}

\par We straightforwardly observe that under the hypotheses of
Proposition \ref{queen3}:

\begin{equation}
\label{lyapunov} (0=) Im^{1|(w,w\cdot\kappa
-K)}_\phi[L_{g}(f)]=\sum_{h\in H_{\sigma+1}} a_h
C^h_{g}(f,\phi)+\sum_{h\in H_{\ge\sigma+2}} a_h C^h_{g}(f,\phi)
\end{equation}
Here the complete contractions indexed in $H_{\sigma+1}$ have
$\sigma+1$ factors while the ones indexed in $H_{\ge\sigma+2}$
have at least $\sigma +2$ factors. All contractions on the RHS are
understood to be in the
 form (\ref{skeleton}); this can be done by decomposing the Weyl tensor.
 We observe that by virtue of the transformation
laws (\ref{weyltrans}), (\ref{levicivita}) and (\ref{wschtrans})
and by virtue of the fact
 that the complete contractions in $L_{g}(f)$ have weight $-K$, all complete
contractions in (\ref{lyapunov}) will have a factor
$\nabla^{(p)}\phi$ with $p\ge 1$. Now, we denote by
$H^{*}_{\sigma+1}\subset H_{\sigma+1}$ the index set of complete
contractions with a factor $\nabla\phi$ (with only one
derivative). We observe, by virtue of the transformation laws
(\ref{weyltrans}), (\ref{levicivita}) and (\ref{wschtrans}), that
this sublinear combination can {\it only} arise from the complete
contractions of length $\sigma$ in $L_{g}(f)$ by applying the
transformation law (\ref{levicivita}) {\it or} by bringing out a
factor $\nabla\phi$ by virtue of the transformations
$W_{ijkl}\rightarrow e^{2\phi}W_{ijkl}$ and $f\rightarrow
e^{w\phi}f$ (when $w\ne 0$).

\par More is true:
Consider $Im^{1|(w,w\cdot\kappa -K)}_\phi[C^u_{g}(f)]$, where
$u\in U^\mu_\sigma$. We pick out the sublinear combination of
complete contractions ({\it in the form (\ref{skeleton})}) with a
factor $\nabla\phi$ and with $\mu-1$ internal contractions in
total. For each $u\in U^\mu_\sigma$ we denote this sublinear
combination by $Im^{1|(w,w\cdot\kappa -K)|*}_\phi[C^u_{g}(f)]$.

\par It is straightforward to observe that for each $u\in U^\mu_\sigma$,
 any complete contraction
in $Im^{1|(w,w\cdot\kappa -K)}_\phi[C^u_{g}(f)]$ (which is written
as a linear combination of complete contractions in the form
(\ref{skeleton})) that does not belong to the sublinear
combination $Im^{1|(w,w\cdot\kappa -K)|*}_\phi[C^u_{g}(f)]$ must
either have a factor $\nabla^{(p)}\phi$ with $p>1$ or must have at
least $\mu$ internal contractions.

\par Furthermore, it is equally straightforward that for any
complete contraction $C^u_g(f)$, $u\in U_\sigma\setminus
U^\mu_\sigma$,
 ($C^u_g(f)$ by hypothesis is in the form (\ref{onlyW}) with $\sigma$
factors and at least $\mu+1$ internal contractions),
$Im^{1|(w,w\cdot\kappa -K)}_\phi[C^u_{g}(f)]$ must consist of
complete contractions in the form (\ref{skeleton}) with either a
factor $\nabla^{(p)}\phi$, $p>1$, or with at least $\mu $ internal
contractions.

\par Therefore, by the above discussion,
we can re-express (\ref{lyapunov}) to obtain:

\begin{equation}
\label{lyapunov2} \sum_{u\in U^\mu_\sigma} a_u
Im^{1|(w,w\cdot\kappa -K)|*}_\phi[C^u_{g}(f)]+\sum_{v\in V_1} a_v
C^v_{g}(f,\phi)+\sum_{v\in V_2} a_v C^v_{g}(f,\phi)=0
\end{equation}
where the complete contractions indexed in $V_1$ are in the
 form (\ref{skeleton}) and have $\sigma+1$
factors and a factor $\nabla^{(p)}\phi$, and moreover either have
$p>1$ or have $p=1$ and at least $\mu$ internal contractions (we
accordingly divide $V_1$ into $V_1^\alpha, V_1^\beta$). The
complete contractions indexed in $V_2$ are generic in the form
(\ref{skeleton}) and have at least $\sigma+2$ factors. Now, since
the above holds at any point $x_0\in M$ and for any values we
assign to the jet of $\phi$ at $x_0$, we derive that:

\begin{equation}
\label{lyapunov2} \sum_{u\in U^\mu_\sigma} a_u
Im^{1|(w,w\cdot\kappa -K)|*}_\phi[C^u_{g}(f)]+\sum_{v\in V^\beta_1} a_v
C^v_{g}(f,\phi)+\sum_{v\in V'_2} a_v C^v_{g}(f,\phi)=0
\end{equation}
(where the contractions indexed in $V'_2$ have at least $\sigma+2$ factors).

 By  Proposition \ref{equiv}, the above must hold formally for
the {\it linearizations} of the complete contractions with
$\sigma+1$ factors (the degree of the complete contractions will
be at most $n+1$--see the note after definition \ref{kakom}). In
the notation of Proposition \ref{equiv}, we will have:

\begin{equation}
\label{lyapunov3} \sum_{u\in U^\mu_\sigma} a_u lin\{
Im^{1|(w,w\cdot\kappa -K)|*}_\phi[C^u_{g}(f)]\}+\sum_{v\in V^\beta_1}
a_v lin\{ C^v_{g}(f,\phi)\} =0
\end{equation}

\par Now, we only have to observe that the
total number of internal contractions in each of the linearized
complete contractions above remain invariant under the
permutations by which we make the above formally zero. Therefore,
(\ref{lyapunov3}) implies:

\begin{equation}
\label{lyapunov4} \sum_{u\in U^\mu_\sigma} a_u lin\{
Im^{1|(w,w\cdot\kappa -K)|*}_\phi[C^u_{g}(f)]\}=0
\end{equation}
and moreover the above still holds formally.

\par We will use the above equation to prove our Proposition \ref{finn},
but in order to do this we must better understand how the
sublinear combination $Im^{1|(w,w\cdot\kappa
-K)|*}_\phi[C^u_{g}(f)]$ arises from $C^u_{g}(f)$, $u\in
U^\mu_\sigma$.
\newline

\par In order to understand $Im^{1|(w,w\cdot\kappa
-K)|*}_\phi[C^u_{g}(f)]$, for any $u\in U^\mu_\sigma$, we will
first seek to separately understand the transformation of each
factor $T$ in $C^u_{g}(f)$ under the operation
$Im^{1|(w,w\cdot\kappa -K)}_\phi[\dots]$. If $T=\nabla^{(m)}W$, we
define $Im^1_\phi[T]$ to stand for the terms of homogeneity one in
$\phi$ in $\nabla^{(m)}W(e^{2\phi}g)$, which have exactly one
derivative on $\phi$; the terms in $Im^1_\phi[T]$ will be in the
form $\nabla^{(m-1)}W\otimes\nabla\phi$ or
$\nabla^{(m-1)}W\otimes\nabla\phi\otimes g$, i.e. we are not
decomposing the Weyl tensor. If $T_y=\nabla^{(p)}f$, we define
$Im^1_\phi[T_y]$ to stand for the terms of homogeneity one in
$\phi$ in $[\nabla^{(p)}(e^{w\cdot\phi}f)](e^{2\phi}g)$, which
have exactly one derivative on $\phi$.

\par It follows that the sublinear combination of terms in
$Im^{1|(w,w\cdot\kappa -K)}_\phi[C^u_{g}(f)]$ with a factor
$\nabla\phi$ arises by picking out each factor $T$ in $C^u_g(f)$,
replacing it by $Im^1_\phi[T]$ and then summing over all these
possible replacements. Now, consider any factor $T_y$ in
$C^u_g(f)$ which has $\delta>0$ internal contraction in
$C^u_g(f)$; define $R[T_y]$ to stand for the sublinear combination
of terms in $Im^1_\phi[T_y]$ {\it in the form $\nabla^{(m-1)}W$ or
$\nabla^{(p-1)}f$} with exactly $\delta-1$ internal contractions.
Define $Im^{1|(w,w\cdot\kappa -K)|\sharp}_\phi[C^u_{g}(f)]$ to
stand for the sublinear combination in $Im^{1|(w,w\cdot\kappa
-K)}_\phi[C^u_{g}(f)]$ which arises by picking out each factor
$T_y$ in $C^u_g(f)$, replacing it by $R[T_y]$ and then summing
over all these substitutions; then observe that
$Im^{1|(w,w\cdot\kappa -K)|*}_\phi[C^u_{g}(f)]$ will be a {\it
sublinear combination} of $Im^{1|(w,w\cdot\kappa
-K)|\sharp}_\phi[C^u_{g}(f)]$. In view of this, we set out to
understand the form of $R[T_y]$ for each factor $T_y$ in
$C^u_g(f)$.
\newline

\par We will introduce some notation to do this. Consider any
$C^u_{g}(f)$, $u\in U^\mu_\sigma$ in the form (\ref{onlyW}). We
pick out the factors $T_y$ that contain an internal contraction,
and assume they are indexed in the set $\{T_1,\dots
,T_{x_u},T_{x_u+1},\dots ,T_{x_u+b_u}\}$. We make the
 convention that the factors $\{T_1,\dots,T_{x_u}\}$ are in the
 form $\nabla^{(m)}W$, while the factors
$T_{x_u+1},\dots T_{x_u+b_u}$ are in the form $\nabla^{(p)}f$.

\par For convenience, we will repeatedly apply the Ricci
 identity to write each complete contraction $C^u_{g}(f)$, $u\in U_\sigma$ in
 the form:

\begin{equation}
\label{onlyW2}
\begin{split}
&contr(\nabla^{y_1\dots y_w}\nabla^{(m_1)}_{r_1\dots r_m}W_{ijkl}
\otimes\dots\otimes\nabla^{h_1\dots h_b}
\nabla^{(m_{\sigma-\kappa})}_{c_1\dots c_{m_{\sigma-\kappa}}}W_{i'j'k'l'}\otimes
\\&\nabla^{f_1\dots f_z}\nabla^{(p_1)}_{q_1\dots q_{p_1}}f
\otimes\dots\otimes \nabla^{d_1\dots d_x}
\nabla^{(p_\kappa)}_{i_1\dots i_{p_\kappa}}f)
\end{split}
\end{equation}
where we are making the convention that the
 raised indices are contracting against one of the lower
 indices in the same factor, and all the lower indices that are not contracting against a
 raised index in the same factor are contracting against some
lower index in another factor. By Lemma \ref{appendix}, we observe
that the correction terms we obtain will satisfy the extra
restrictions when they are relevant.

We first consider any factor $T_y, y\in \{x_u+1,\dots x_u+b_u\}$
that has, say, $\delta$ internal contractions. Recall $T_y$ in the
form $T_y=\nabla^{u_1\dots u_\delta} \nabla^{(p)}_{r_1\dots r_p}f$
with the conventions that each of
 the indices ${}^{u_1},\dots ,{}^{u_\delta}$
is contracting against
one of the indices ${}_{r_1},\dots ,{}_{r_p}$. Moreover each
 of the indices ${}_{r_1},\dots ,{}_{r_p}$ that is not
 contracting
  against an index ${}^{u_1},\dots ,{}^{u_\delta}$ is a
 free index. It follows that (modulo terms with three factors) $R[T_y]$ will be:

 \begin{equation}
 \label{mert1}
R[T_y]= [\delta\cdot (n-2)-4{\delta\choose{2}}+2\delta w]
(\nabla)^{u_1}\phi\nabla^{u_2\dots u_\delta} \nabla^{(p)}_{r_1\dots
r_p}f \end{equation}

\par Now, we pick out any factor $T_y,y\in \{1,\dots x_u\}$, in the form $\nabla^{u_1\dots u_\delta}
\nabla^{(m)}_{r_1\dots r_m}W_{ijkl}$. In order to understand
$R[T_y]$ in this setting we must distinguish three cases: In the
 first case we have that all the internal contractions in
$T_y$ are between derivative indices. In the second case we have
that there is precisely one internal contraction involving an
internal index (and with no loss of generality we will assume that
${}^{u_1}$ is contracting against the index ${}_i$ in $W_{ijkl}$).
In the third case we have there are precisely two internal
contractions
 involving internal indices, (and with no loss of
  generality we will assume that the indices
${}^{u_1},{}^{u_2}$ are contracting against the indices
${}_i,{}_k$ respectively).

 In the case where all the internal
contractions in $T_y$ are between derivative indices we see that
(modulo terms with 3 factors) $R[T_y]$ will be:

\begin{equation}
\label{mert2} R[T_y]= [\delta\cdot (n-2)-4{\delta\choose{2}}]
\nabla^{u_2\dots u_\delta}\nabla^{(m)}_{r_1\dots r_m}W_{ijkl}
(\nabla)^{u_1}\phi \end{equation}

\par In the case where there is precisely one internal contraction
in $T_y$ (between the indices ${}^{u_1},{}_i$), we also assume for
convenience that if $\delta>1$ then ${}^{u_2}$ is contracting
against the index ${}_{r_1}$. Then, (modulo terms with 3 factors)
$R[T_y]$ will be the sum:

\begin{equation}
\label{erdogan}
\begin{split}
&(n-3)\nabla^{u_2\dots u_\delta} \nabla^{(m)}_{r_1\dots
r_m}W_{ijkl}(\nabla)^i\phi+ [(\delta-1)(n-2)-4{\delta\choose{2}}]
\nabla^{u_1u_3\dots u_\delta}\nabla^{(m)}_{r_1\dots r_m}
W_{ijkl}(\nabla)^{u_2}\phi+
\\&(\delta-1)\nabla^{u_1u_3\dots u_\delta}
\nabla^{(m)}_{jr_2\dots r_m} W_{ir_1kl}(\nabla)^{r_1}\phi
\end{split}
\end{equation}

\par Now, we consider the case where a factor $T_i$ is in the
 form $\nabla^{u_1\dots u_\delta}\nabla^{(m)}_{r_1\dots r_m}W_{ijkl}$ where two
 of the internal indices ($i$ and $k$) are involved in an
 internal contraction. We then assume for convenience that if $\delta>2$ then
  ${}^{u_3}$ is
contracting against ${}_{r_1}$. Thus, modulo terms with 3 factors, $R[T_i]$ will be:

\begin{equation}
\label{pareme}
\begin{split}
&(n-4)\nabla^{u_2\dots u_\delta} \nabla^{(m)}_{r_1\dots
r_m}W_{ijkl}(\nabla)^i\phi+ (n-4)\nabla^{u_1u_3\dots u_\delta}
\nabla^{(m)}_{r_1\dots r_m}W_{ijkl}(\nabla)^k\phi+
\\& (\delta-2)(n-2)\nabla^{u_1u_2u_4\dots u_x}
\nabla^{(m)}_{r_1\dots r_m}W_{ijkl}(\nabla)^{r_1}\phi
\\&-4[{{\delta-2}\choose{2}}+2(\delta-2)]
\nabla^{u_1u_2u_4\dots u_x} \nabla^{(m)}_{r_1\dots
r_m}W_{ijkl}(\nabla)^{r_1}\phi+
\\&(\delta-2)\nabla^{u_1u_2u_4\dots u_\delta}
\nabla^{(m)}_{jr_2\dots r_m} W_{ir_1kl}(\nabla)^{r_1}\phi+
(\delta-2)\nabla^{u_1u_2u_4\dots u_\delta} \nabla^{(m)}_{lr_2\dots
r_m} W_{ijkr_1}(\nabla)^{r_1}\phi
\end{split}
\end{equation}

\par Now, for each $y\in \{1,\dots x_u+b_u\}$
we define $C^{u,y}_{g}(f,\phi)$ to be the complete contraction (or
linear combination of complete contractions), in the form
(\ref{skeleton}), that arises from $C^u_{g}(f)$ by
 replacing the factor $T_y$ by $R[T_y]$ and then replacing each factor
$\nabla^{(m)}_{r_1\dots r_m}W_{ijkl}$ by a factor
$\nabla^{(m)}_{r_1\dots r_m}R_{ijkl}$ or $\frac{n-3}{n-2}
\nabla^{(m)}_{r_1\dots r_m}R_{ijkl}$ depending on whether there is
no internal contraction involving an internal index or not.

\par Then, for each $u\in U^\mu_\sigma$ we define
$C^{u,+}_{g}(f,\phi)$ to stand for the linear combination of
complete contractions:

\begin{equation}
\label{gosta} C^{u,+}_{g}(f,\phi)=\sum_{y=1}^{x_u+b_u}
C^{u,y}_{g}(f,\phi)
\end{equation}

\par Now, by the definitions of $R[T_y]$ and
$Im^{1|(w,w\cdot\kappa -K)|*}_\phi[C^u_{g}(f)]$, by the formulas
above and by the decomposition formulas for $\nabla^{(m)}W$, we
derive that if $C^u_{g}(f)$, $u\in U^\mu_\sigma$, has no factor
$\nabla^{(m)}W$ with two internal contractions involving internal
indices then:

\begin{equation}
\label{leffler} Im^{1|(w,w\cdot\kappa
-K)|*}_\phi[C^u_{g}(f)]=C^{u,+}_{g}(f,\phi)
\end{equation}
On the other hand, if $C^u_{g}(f)$ does contain factors
$\nabla^{(m)}W$ with two internal contractions involving internal
indices then:

\begin{equation}
\label{leffler2} Im^{1|(w,w\cdot\kappa
-K)|*}_\phi[C^u_{g}(f)]=C^{u,+}_{g}(f,\phi)+ \sum_{d\in D} a_d
C^d_{g}(f,\phi)
\end{equation}
where $\sum_{d\in D} a_d C^d_{g}(f)$ stands for a generic linear
combination of complete contractions in the form (\ref{skeleton})
with the extra feature that they contain $a>0$ factors
$\nabla^z_{r_1\dots r_z}S$, where $S$ is the
 scalar curvature. The contractions $C^d_{g}(f,\phi)$ arise by
 virtue of the factors $T^z, z\in Z^{\delta=x}$ in
(\ref{decompo2}).

\par Thus, equation (\ref{lyapunov4}) implies:

\begin{equation}
\label{lyapunov5} \sum_{u\in U^\mu_\sigma}a_u lin\{
C^{u,+}_{g}(f,\phi)\}+\sum_{d\in D} a_d lin\{ C^d_{g}(f,\phi)\}=0
\end{equation}

\par We then have two claims.

\begin{lemma}
\label{simp1} Given (\ref{lyapunov5}) we claim that:
\begin{equation}
\label{lyapunov6} \sum_{d\in D} a_d lin\{ C^d_{g}(f,\phi)\}=0
\end{equation}
\end{lemma}

\par Observe that if we can show the above, we will then have:

\begin{equation}
\label{lyapunov7} \sum_{u\in U^\mu_\sigma}a_u lin\{
C^{u,+}_{g}(f,\phi)\}=0
\end{equation}
We present our next claim with a little more notation: We denote
by $lin\{C^{u,+,i_1\dots i_{\mu-1}}_{g}(f,\phi)\}$ the
(linearized) tensor field that arises from $lin\{
C^{u,+}_{g}(f,\phi)\}$ by replacing each of the internal
contractions by a free index (meaning that in each complete
contraction $(\nabla^a,{}_a)$
 we erase $\nabla^a$ and make the index
${}_a$ free). We then form a
 (linearized) complete contraction
$lin\{C^{u,+,i_1\dots i_{\mu-1}}_{g}(f,\phi)\nabla_{i_1}\upsilon
\dots\nabla_{i_{\mu-1}}\upsilon\}$. Our second claim is that:

\begin{lemma}
\label{simp2} Assuming (\ref{lyapunov7}) and employing the
notation above we have that:
\begin{equation}
\label{lyapunov8} \sum_{u\in U^\mu_\sigma}a_u lin\{C^{u,+,i_1\dots
i_{\mu-1}}_{g}(f,\phi) \nabla_{i_1}\upsilon
\dots\nabla_{i_{\mu-1}}\upsilon\}=0
\end{equation}
and moreover the above holds formally.
\end{lemma}

\par We will show these two Lemmas below. For now, let us
 check how these Lemmas will imply  Proposition \ref{finn}
 (and hence  Proposition \ref{queen3}). We arbitrarily pick out an element
$\vec{\alpha}\in \Lambda$ (see definition \ref{tereza}) and we
will show Proposition \ref{finn} for the index
 set $U^{\mu,\vec{\alpha}}_\sigma$. We must distinguish two
 cases: Either for the character $\vec{\alpha}=(RL_1|RL_2)$
  we have $RL_1=\emptyset$
(in which case we must also necessarily have $RL_2\ne\emptyset$,
since $\mu>0$), or we have $RL_1\ne\emptyset$. For future
reference, we will write out $\vec{\alpha}$ explicitly:

\begin{equation}
\label{indos} \vec{\alpha}=(RL_1|RL_2)=(\zeta_1,\dots ,\zeta_a
|\xi_1,\dots ,\xi_b)
\end{equation}

\par Now, we first consider the case $RL_1=\emptyset$ for our chosen
$\vec{\alpha}\in \Lambda$. As we observed, this implies that
$RL_2\ne\emptyset$, i.e. for each $u\in
U^{\mu,\vec{\alpha}}_\sigma$ we have $\xi_1>0$.

 We then
consider the sublinear combination in (\ref{lyapunov8}) which
contains complete contractions with the (linearized) factor
$\nabla\phi$ contracting against a (linearized) factor $T=lin\{
\nabla^{(p)}f\}$, and moreover the factor $T$ is also contracting
against precisely $\xi_1-1$ factors $\nabla\upsilon$. If we denote
the linear combination of such complete contractions by $lin\{
S_{g}(f,\phi,\upsilon)\}$ then since (\ref{lyapunov8}) holds
formally, we derive that:

\begin{equation}
\label{piskopio} lin\{S_{g}(f,\phi,\upsilon)\}=0
\end{equation}
and moreover the above holds formally. Now, we set $\phi=\upsilon$
in the above. By applying equations (\ref{lyapunov4}) and
(\ref{gosta}) we observe that:

\begin{equation}
\label{lyapunov12} \begin{split} &
(0=)lin\{S_{g}(f,\upsilon,\upsilon)\}= Const(\vec{\alpha})
\sum_{u\in U^{\mu,\vec{\alpha}}_\sigma}a_u lin\{C^{u|i_1\dots
i_{\mu}}_{g}(f) \nabla_{i_1}\upsilon
\dots\nabla_{i_{\mu}}\upsilon\}
\\&+\sum_{s\in S} a_s lin\{ C^{s,i_1\dots
i_\mu}_{g}\nabla_{i_1}\upsilon\dots \nabla_{i_\mu}\upsilon\}
\end{split}
\end{equation}
($C^{u|i_1\dots i_{\mu}}_{g}$ stands for the tensor field
 that arises from $C^u_{g}(f)$ by making all complete
 contractions into free indices). Here each $lin\{
C^{s,i_1\dots i_\mu}_{g}\nabla_{i_1}\upsilon\dots
\nabla_{i_\mu}\upsilon\}$ has a character that is {\it different}
from our chosen $\vec{\alpha}$. Therefore, since the above holds
formally and since the character of the complete contractions
remains unaltered under the permutations that make
(\ref{lyapunov12}) {\it formally} zero,  we derive:

\begin{equation}
\label{kutlu3} Const(\vec{\alpha}) \sum_{u\in
U^{\mu,\vec{\alpha}}_\sigma}a_u lin\{C^{u|i_1\dots i_{\mu}}_{g}(f)
\nabla_{i_1}\upsilon \dots\nabla_{i_{\mu}}\upsilon\}=0
\end{equation}
The constant $Const(\vec{\alpha})$  only depends on our chosen
character $\vec{\alpha}$. To be precise, if we denote by $m$ the
number of times that $\xi_1$ appears in $\vec{\alpha}$ ($m\ge 1$),
by applying (\ref{mert1}) and the definitions, it follows that:

$$Const(\vec{\alpha})=m\cdot[\xi_1\cdot(n-2)-4{{\xi_1}
\choose{2}}+2\xi_1\cdot w]=m\cdot \xi_1\cdot[n-2\xi_1+2w]$$ and we
observe that $Const'(\vec{\alpha})\ne 0$ if $(w+\frac{n}{2})\notin
\mathbb{Z}_{+}$. Furthermore, if $w=-\frac{n}{2}+k$ for some $k\in
\mathbb{Z}_{+}$, we will still have $Const(\vec{\alpha})\ne 0$ by
our extra restriction on $\beta[L_g(f)]$, which ensures that
$\xi_1<k$.
\newline

\par Now, we proceed with the second case, where $RL_1\ne \emptyset$.
Thus, for every $u\in U^{\mu,\vec{\alpha}}_\sigma$ there is an
internal contraction in some factor $\nabla^{(m)}W$. Now,
similarly to the case above, we will consider the sublinear
combination of contractions in (\ref{lyapunov8}) for which
$\nabla\phi$ is contracting against a linearized curvature factor,
say $T$, and also $T$ is contracting against precisely
$(\zeta_1-1)$ factors $\nabla\upsilon$.

If we denote the linear combination of such complete contractions
by $lin\{ S'_{g}(f,\phi,\upsilon)\}$ then since (\ref{lyapunov8})
holds formally, we derive that:

\begin{equation}
\label{piskopio} lin\{S'_{g}(f,\phi,\upsilon)\}=0
\end{equation}
and moreover the above holds formally. Now, we set $\phi=\upsilon$
and derive a new equation which we denote:
$lin\{S'_{g}(f,\upsilon,\upsilon)\}=0$ (and the above still holds
formally). By the same argument as for the previous case, we
derive that:

\begin{equation}
\label{ramirez}
\begin{split}
& (0=)lin\{S'_{g}(f,\phi,\upsilon)\}= Const'(\vec{\alpha})
\sum_{u\in U^{\mu,\vec{\alpha}}_\sigma}a_u lin\{C^{u|i_1\dots
i_{\mu}}_{g}(f) \nabla_{i_1}\upsilon
\dots\nabla_{i_{\mu}}\upsilon\}
\\&+\sum_{s\in S} a_s lin\{ C^{s,i_1\dots
i_\mu}_{g}\nabla_{i_1}\upsilon\dots \nabla_{i_\mu}\upsilon\}
\end{split}
\end{equation}
(here $C^{u|i_1\dots i_{\mu}}_{g}$ stands for the tensor field
 that arises from $C^u_{g}(f)$ by making all internal
 contractions into free indices). Here again the contractions $lin\{
C^{s,i_1\dots i_\mu}_{g}\nabla_{i_1}\upsilon\dots
\nabla_{i_\mu}\upsilon\}$ have a {\it different} character from
our chosen $\vec{\alpha}$. Therefore, since the above holds
formally, we derive:

\begin{equation}
\label{kutlu3} Const'(\vec{\alpha}) \sum_{u\in
U^{\mu,\vec{\alpha}}_\sigma}a_u lin\{C^{u|i_1\dots i_{\mu}}_{g}(f)
\nabla_{i_1}\upsilon \dots\nabla_{i_{\mu}}\upsilon\}=0
\end{equation}
The constant $Const'(\vec{\alpha})$ only depends on our chosen
character $\vec{\alpha}$. In particular, if we denote by $f$ the
number of times that the number $\zeta_1$ appears in
(\ref{indos}), then it follows that:
\begin{equation}
\label{exupery} Const'(\vec{\alpha})= f\cdot [\zeta_1\cdot
(n-2)-4{{\zeta_1}\choose{2}}]=f\cdot\zeta_1\cdot (n-2\zeta_1)
\end{equation}
 (\ref{exupery}) follows
from (\ref{mert1}), (\ref{mert2}), (\ref{erdogan}),
(\ref{pareme}), (\ref{gosta}). Now, in the case where $n$ is odd,
we clearly see that $Const'(\vec{\alpha})\ne 0$. In $n$ is even,
we have only have to recall that $\zeta_1$ corresponds to the
number of internal contractions in some  factor $T=\nabla^{(m)}W$
in some complete contraction in $L^\sigma_g(f)$. But then our restriction
$\gamma[L_g(f)]<\frac{n}{2}$ guarantees that $\zeta_1<\frac{n}{2}$
(observe that $\gamma[T]\ge \zeta_1$). Thus we derive
$Const'(\vec{\alpha})\ne 0$ in this case also.

\par In view of the above we also derive our claim in the case
where $RL_1\ne\emptyset$ for our chosen $\vec{\alpha}\in \Lambda$.
$\Box$

\par In view of the above, it suffices to show our Lemmas
\ref{simp1} and \ref{simp2} to derive Proposition \ref{finn} (and
hence theorems \ref{ayto}, \ref{ayto2}, \ref{careful}). We will
start by proving Lemma \ref{simp2}, because Lemma \ref{simp1} is
slightly more complicated.
\newline

{\it Proof of Lemma \ref{simp2}:}
\newline

\par In order to see this claim it will be more useful to
 consider complete contractions in the curvature and its
 covariant derivatives, i.e.~ in the form (\ref{skeleton}), rather than contractions in
 linearized tensors.

\par Our point of departure is equation (\ref{lyapunov7}). We
``memorize'' the sequence of permutations that make the left hand
side of (\ref{lyapunov7}) formally zero.  Now, we consider the
complete contractions $C^{u,+}_{g}(f,\phi)$ in any higher
dimension $N$, so we obtain complete contractions
$C^{u,+}_{g^N}(f,\phi)$, $u\in U^\mu_\sigma$. Now, we repeat the
sequence of permutation that made (\ref{lyapunov}) formally zero
to the linear combination $\sum_{u\in U^\mu_\sigma} a_u
C^{u,+}_{g^N}(f,\phi)$ and we derive a new equation:

\begin{equation}
\label{lyapunov5unsymN} \sum_{u\in U^\mu_\sigma}a_u
C^{u,+}_{g^N}(f,\phi)=\sum_{p\in P} a_p C^{p}_{g^N}(f,\phi)
\end{equation}
for any $N\ge n$. Here the complete contractions
$C^{p}_{g^N}(f,\phi)$ are in the
 form (\ref{skeleton}) and have at least $\sigma +2$ factors. The
 contractions $C^{p}_{g^N}(f,\phi)$ arise as correction terms in
(\ref{lyapunov5unsymN}) due to the right hand side of
(\ref{curvature}).
 Therefore, we
derive that each $C^{p}_{g}(f,\phi)$ will not contain a factor $S$
(of the scalar curvature), because such factor cannot arise in
correction terms appearing by virtue of the curvature identity.

\par One more piece of notation. For any complete contraction
$C_{g^N}(f)$ (in any dimension $N$) of weight $-K$ we define the
$l^{th}$ conformal variation
 $Var^l_\upsilon$ ($l\in \mathbb{Z}_{+}$ and $\upsilon$ is an
arbitrary scalar function):

$$Var^l_\upsilon[C_{g^N}(f,\phi)]=\frac{\partial^l}{\partial\lambda^l}|_{\lambda=0}
[C_{e^{2\lambda\upsilon g^N}}(f,\phi)]$$

\par We now consider $Var^{\mu-1}_{\upsilon}$ of (\ref{lyapunov5unsymN}).
 Clearly, we have that:

\begin{equation}
\label{lyapunov5unsym} Var^{\mu-1}_\upsilon \{ \sum_{u\in
U^\mu_\sigma}a_u C^{u,+}_{g^N}(f,\phi)\}=Var^{\mu-1}_\upsilon
\{\sum_{p\in P} a_p C^{p}_{g^N}(f,\phi)\}
\end{equation}
\par Now, by a careful study of the transformation laws
(\ref{levicivita}) and (\ref{weyltrans}) we observe that we can
write:

\begin{equation}
\label{ooo} \begin{split} & Var^{\mu-1}_\upsilon \{ \sum_{u\in
U^\mu_\sigma}a_u C^{u,+}_{g^N}(f,\phi)\}=N^{\mu-1} \sum_{u\in
U^\mu_\sigma}a_u C^{u,+,i_1\dots
i_{\mu-1}}_{g^N}(f,\phi)\nabla_{i_1}\upsilon\dots
\nabla_{i_{\mu-1}}\upsilon+
\\&N^{\mu-1}\sum_{t\in T} a_t
C^t_{g^N}(f,\phi,\upsilon)+\sum_{y=0}^{\mu-2}N^y \sum_{h\in H^y}
a_h C^h_{g^N}(f,\phi,\upsilon)
\end{split}
\end{equation}
All the complete contractions above are in the form:
\begin{equation}
\label{monoagkalia}
\begin{split}
&contr(\nabla^{(m_1)}R\otimes\dots\otimes\nabla^{(m_t)}R\otimes
\\&\nabla^{(\nu_1)}f\otimes\dots\otimes\nabla^{(\nu_\kappa)}f
\otimes\nabla^{(p_1)}\upsilon\otimes\dots\otimes\nabla^{(p_{\mu-1})}\upsilon)
\end{split}
\end{equation}
 here the
complete contractions $C^t_{g}(f,\phi,\upsilon)$ have $\sigma +\mu$
factors, and they are in the form (\ref{monoagkalia})  and
moreover have the property that at least one $p_i,i=1,\dots
,\mu-1$ is $\ge 2$. $\sum_{h\in H^y}\dots$ stands for a generic
linear combination of complete contractions in the above form
(what is important is that it is multiplied by $N^\beta$ with
$\beta<\mu-1$). Furthermore, in the notation above, we also
 have:

\begin{equation}
\label{ooo2}
\begin{split} & Var^{\mu-1}_\upsilon \{ \sum_{p\in
P} a_p C^{p}_{g^N}(f,\phi)\}=
\\&N^{\mu-1}\sum_{t\in T'} a_t
C^t_{g^N}(f,\phi,\upsilon)+\sum_{y=0}^{\mu-2}N^y \sum_{h\in H^y}
a_h C^h_{g^N}(f,\phi,\upsilon)
\end{split}
\end{equation}
where the complete contractions indexed in $T'$ are in the
 form (\ref{monoagkalia}) but with at least $\sigma+\mu+1$
factors.

Therefore, by virtue of (\ref{lyapunov5unsym})
(which holds formally, for $N$ large enough) and the analysis
 above, we obtain:
\begin{equation}
\label{ooo3} \begin{split} & N^{\mu-1} \sum_{u\in U^\mu_\sigma}a_u
C^{u,+,i_1\dots i_{\mu-1}}_{g^N}(f,\phi)\nabla_{i_1}\upsilon\dots
\nabla_{i_{\mu-1}}\upsilon+
\\&N^{\mu-1}\sum_{t\in T\bigcup T'} a_t
C^t_{g^N}(f,\phi,\upsilon)+\sum_{y=0}^{\mu-2}N^y \sum_{h\in H^y}
a_h C^h_{g^N}(f,\phi,\upsilon)=0
\end{split}
\end{equation}
and furthermore this equation holds formally for $N\ge n+\mu$.
Therefore, by just picking $M^N=M^{n+\mu}\times
S^1\times\dots\times S^1$ we obtain:

\begin{equation}
\label{ooo3} \begin{split} & N^{\mu-1} \sum_{u\in
U^\mu_\sigma}a_u C^{u,+,i_1\dots
i_{\mu-1}}_{g^{n+\mu}}(f,\phi)\nabla_{i_1}\upsilon\dots
\nabla_{i_{\mu-1}}\upsilon+
\\&N^{\mu-1}\sum_{t\in T\bigcup T'} a_t
C^t_{g^{n+\mu}}(f,\phi,\upsilon)+\sum_{y=0}^{\mu-2}N^y
\sum_{h\in H^y} a_h C^h_{g^{n+\mu}}(f,\phi,\upsilon)=0
\end{split}
\end{equation}
and the above still holds formally. Here $N$ can be any integer
with $N\ge n+\mu$. Thus, if we treat
 the above as a polynomial in $N$, we derive that:

\begin{equation}
\label{ooo4} \begin{split} &  \sum_{u\in U^\mu_\sigma}a_u
C^{u,+,i_1\dots
i_{\mu-1}}_{g^{n+\mu}}(f,\phi)\nabla_{i_1}\upsilon\dots
\nabla_{i_{\mu-1}}\upsilon+ \sum_{t\in T\bigcup T'} a_t
C^t_{g^{n+\mu}}(f,\phi,\upsilon)=0
\end{split}
\end{equation}
and the above holds formally. Therefore, it must also hold
formally at the linearized level, so we derive:

\begin{equation}
\label{ooo5} \begin{split} &  \sum_{u\in U^\mu_\sigma}a_u lin\{
C^{u,+,i_1\dots i_{\mu-1}}_{g}(f,\phi)\nabla_{i_1}\upsilon\dots
\nabla_{i_{\mu-1}}\upsilon\}=0
\end{split}
\end{equation}
This is precisely our desired conclusion. $\Box$
\newline

{\it Proof of Lemma \ref{simp1}:} We divide the index set $D$ into
subsets $D^a, a=1,\dots ,\sigma -\kappa$, according to the number
of factors $\nabla^{(a)}S$ that the contractions $C^d_{g}(f), d\in
D$ contains. It would clearly suffice to show that for each
$c=1,\dots ,\sigma-\kappa$ we must have:

\begin{equation}
\label{mine} \sum_{d\in D^c} a_d lin\{ C^d_{g}(f,\phi)\}=0
\end{equation}

\par We show (\ref{mine}) by a downward induction on $c$:
 We assume that for some number $M\le \sigma-\kappa$ we have that
for every $c>M$ $\sum_{d\in D^c} a_d lin\{ C^d_{g}(f)\}=0$. We
 will then show that this must be true for $c=M$.

\par In order to see this claim it will again be more useful to
consider the {\it non-linearized} version of (\ref{lyapunov5}); in
any dimension $N\ge n$:

\begin{equation}
\label{lyapunov5unsym'N} E_{g^N}(f)= \sum_{u\in U^\mu_\sigma}a_u
C^{u,+}_{g^N}(f,\phi)+\sum_{d\in D} a_d  C^d_{g^N}(f,\phi)-
\sum_{p\in P} a_p C^{p}_{g^N}(f,\phi)=0
\end{equation}
where the complete contractions $C^{p}_{g^N}(f,\phi)$ are in the
 form (\ref{skeleton}) and have at least $\sigma +2$ factors.

\par Now, we take the $Var^{M}_{Y}[E_{g^N}(f,\phi)]$
in the above equation, so we will have:

\begin{equation}
\label{lebentes} Im^{M}_{Y}[E_{g^N}(f,\phi)]=0
\end{equation}
We want to describe the left hand side of the above. For each
$d\in D^M$ let us
 denote by $C^d_{g^N}(f,\phi,Y)$
the tensor field that arises from $C^d_{g^N}(f,\phi)$ by replacing
 each of the $M$ factors $\nabla^{(p)}_{r_1\dots r_p}S$ by a
 factor $-2\nabla^{(p)}_{r_1\dots r_p}\Delta Y$. By virtue of the
 transformation laws (\ref{curvtrans}) and
 (\ref{levicivita}), and by applying the same argument as for the previous case we derive that:

\begin{equation}
\label{grafe} \sum_{d\in D^M} a_d  C^d_{g^N}(f,\phi,Y)=0
\end{equation}
modulo complete contractions of greater length, and moreover the
above holds formally.

 We then
 denote by $C^{d,i_1\dots i_{\mu-M-1}}_{g^N}(f,\phi,Y)$ the
 tensor field
 that arises by making each of the internal contractions in
$C^d_{g}(f,\phi,Y)$ into a free index.

\par By the same argument as for the previous Lemma, we derive an
 equation:
\begin{equation}
\label{grafe2} \sum_{d\in D^M} a_d  lin\{ C^{d,i_1\dots
i_{\mu-1-M}}_{g^N}
(f,\phi,Y)\nabla_{i_1}\upsilon\dots\nabla_{i_{\mu-M-1}}
\upsilon\}=0
\end{equation}
 and the above holds formally.
 Now, we observe that in each factor $\nabla^{(A)}Y$ ($A\ge 1$), the last index is
 contracting against a factor $\nabla\upsilon$. Thus, it is
 not difficult to derive that we can make the above linear
 combination formally zero {\it without} permuting the last
 index in each of the $M$ factors
$\nabla^{(p)}Y$. Now, we define an operation $Op\{\dots\}$ that acts
on the complete contractions above by replacing
 each expression
$\nabla^{(p)}_{r_1\dots r_p}Y(\nabla)^{r_p}\upsilon$ (i.e. we pick
out a factor $\nabla^{(p)}_{r_1\dots r_p}Y$ and also the
 factor $\nabla\upsilon$ that is contracting against its last
 index) by an expression $\frac{1}{2}\nabla^{(p-1)}_{r_1\dots r_{p-1}}S$.
Furthermore, we replace each remaining factor $\nabla\upsilon$ by
an internal contraction--i.e. we erase each factor
$\nabla_r\upsilon$ that is contracting against some factor $T$
($T=\nabla^{(p)}f$ or $T=\nabla^{(m)}R$ or $T=\nabla^{(y)}S$) and
then we replace $T$ by $\nabla^rT$.
 Then, since (\ref{grafe2}) holds formally and without permuting
the last index $r_P$ in each
 factor $\nabla^{(p)}_{r_1\dots r_p}Y$, just by repeating these
permutations we derive that:

\begin{equation}
\label{grafe3} \sum_{d\in D^{m}} a_d  Op[lin\{
 C^{d,i_1\dots i_{\mu-1-M}}_{g}(f,Y)\nabla_{i_1}
 \upsilon\dots\nabla_{i_{\mu-1-M}}\upsilon\}]=0
\end{equation}
This is precisely our desired conclusion. $\Box$

\section{Sketch of the proof of Theorem \ref{easy}.}

\par As mentioned, this proof closely follows the methods
developed in \cite{beg:itccg}. All facts that we claim without
proof come from \cite{beg:itccg}. We briefly present the argument
that appears in \cite{a:thesis}.

\par The first step is to express the Riemannian operator
$L_{g}(f)$ as a polynomial in the components of the ambient
curvature $\tilde{R}$
 and its covariant
 derivatives $\tilde{\nabla}^{(m)}\tilde{R}$, and in the harmonic extension $\tilde{u}_w$
  and its covariant derivatives $\tilde{\nabla}^{(p)}\tilde{u}_w$. This can be done by using Graham's
   {\it conformal normal scale}: We have that for any $x_0\in M$ we can pick a metric
   $g_1\in [g]$ so that $\nabla^{(p)}_{(r_1\dots r_p}Ric_{r_{p+1}r_{p+2})}(x_0)=0$ for $p\le N$,
    where we are free to pick $N$ as large as we like.

 \par Now, since $L_{g}(f)$ is assumed to be conformally invariant, it would
 suffice to show that $L_{g_1}(f)$ can be written as a Weyl operator (the conformal invariance will then
 imply that $L_{g}(f)$ can also be written as a Weyl operator).
Therefore we consider $L_{g_1}(f)$ and perform the ambient metric
construction for $(M,g_1)$, as explained in subsection \ref{amb}:
We locally embed the manifold $(M,g_1)$ into the ambient
pseudo-Riemannian $(\tilde{G}_1,\tilde{g}_1)$, mapping the point
$x_0\in M$ to $x_{*}=(1,x_0,0)\in \tilde{G}$, and we consider the
extension of the density $f_w$ to a $w$-homogeneous harmonic
function $\tilde{u}_w$ on $(\tilde{G}_1,\tilde{g}_1)$. Recall that
since $(w+\frac{n}{2})\notin \mathbb{Z}_{+}$,
 this extension is well-defined to any order.  Therefore, we obtain that there exists a
 fixed polynomial $\Pi(\{ \tilde{\nabla}^{(m)}\tilde{R}\}, \{\tilde{\nabla}^{(p)}\tilde{u}_w\})$
  in the components of the tensors $\tilde{\nabla}^{(m)}\tilde{R}$ and $\tilde{\nabla}^{(p)}\tilde{u}_w$ so that:

\begin{equation}
\label{flynt} L_{g_1}(f)=\Pi(\{ \tilde{\nabla}^{(m)}\tilde{R}\},
\{\tilde{\nabla}^{(p)}\tilde{u}_w\})
\end{equation}

\par Now, consider any conformal transformation
$\psi:(M,g_1)\rightarrow (M,g_2)$ with $\psi(x_0)=x_0$,
$\nabla\psi(x_0)=m^i_j\in O(g(x_0))$, for which
$\psi^{*}g_2=e^{2\phi}g_1$. Now, perform the ambient metric
construction $(\tilde{G}_2,\tilde{g}_2)$ for $g_2$, mapping $x_0$
to the point $x_{*}=(1,x_0,0)\in\tilde{G}_2$. As discussed in
subsection \ref{amb}, there then exists an isometry (mod
$O(\rho^\infty)$), $\Phi:\tilde{G}_2\rightarrow \tilde{G}_1$, for
which $\Phi(1,x_0,0)=(e^{x_0},x_0,0)$ and moreover $\nabla
\Phi(1,x_0,0)$ is given by the matrix:

\begin{equation}
\label{transPhi}
\left[
\begin{array}{ccc}
        \lambda &\omega_i&t\\
        0       &m^i_j   & s^i\\
        0       &0       &\lambda^{-1}
\end{array}
\right]
\end{equation}
where we set $\lambda(x)=e^{\phi(x)}$ and then $\lambda
=\lambda(x_0)$, $\omega_i=\nabla_i\lambda(x_{0})$,
  $m^i_j\in O(g(x_0))$,
$t=-\frac{1}{2\lambda}\omega_j \omega^j$ (note that $t$ is {\it
not} the coordinate function from above-we are just following the
notations from \cite{beg:itccg}),
$s^i=-\frac{1}{\lambda}m^{ij}\omega_j$. Denote this matrix by $A$.
 Thus, under a conformal
  re-scaling above, we  have that the components
   of the tensors $\tilde{\nabla}^{(m)}_{a_1\dots a_m}\tilde{R}_{a_{m+1}\dots a_{m+4}}$
    and $\tilde{\nabla}^{(p)}_{b_1\dots b_p}\tilde{u}_w$ for
    $\tilde{g}_2$ at $(1,x_0,0)$ arise from the components of
    those tensors of $\tilde{g}_1$ at $(e^{\phi(x_0)}, x_0,0)$ by
  multiplying the vectors $X_i$  by the matrix $A$.
Denote the polynomial that is thus obtained by $\Pi(\{ A\cdot
\tilde{\nabla}^{(m)}\tilde{R}\}, \{A\cdot
\tilde{\nabla}^{(p)}\tilde{u}_w\})$.

  \par On the other hand, the conformal invariance of $L_{g_1}(f)$ guarantees that:

  \begin{equation}
  \label{garanti}
  \begin{split}
 & \Pi(\{A\cdot \tilde{\nabla}^{(m)}\tilde{R}\},
 \{A\cdot \tilde{\nabla}^{(p)}\tilde{u}_w\})(e^{\phi(x_0)},x_0,0)=
 [L_{g_2}(e^{w\phi(x)}f)](x_0)=
 \\&  [e^{(\kappa\cdot w-K)\phi(x)}L_{g_1}(f)](x_0)=
 e^{(\kappa\cdot w-K)\phi(x)}\Pi(\{ \tilde{\nabla}^{(m)}\tilde{R}\}, \{\tilde{\nabla}^{(p)}\tilde{u}_w\})
 \end{split}
 \end{equation}

 \par Thus, the conformal invariance of $L_{g}(f)$ translates into an invariance of the polynomial $\Pi$
 (up to multiplying by a power of $\lambda$) under the action of the
 group $P$ of matrices in the form (\ref{transPhi}). We will call this
 property $P$-invariance.
\newline

\par Now, it is more convenient to work at the linearized level of modules.

  We consider the vector space $W=\mathbb{R}^{n+2}$ and we denote by

\noindent $X^I=( X^0,X^1,\dots ,X^n,X^\infty)$ any point in $W$.
We recall the
 metric $\tilde{g}(x_{*})$ from the subsection \ref{amb}, which
 now defines a quadratic form $\tilde{g}$ on $W$:
 $\sum_{1\le I,J\le n} g_{IJ}X^I
X^J+2X^0X^{\infty}$. We denote by $Q$ the light cone in $W$, with
respect to this quadratic form. We consider the null vector
$e_0=(1,0,\dots ,0)^t\in W$. Now, if we denote by $G$ the
identity-connected component of $O(\tilde{g})$ then the group of
matrices $P$ in the form (\ref{transPhi}) above stands for  the
parabolic subgroup of $G$:

$$P=\{ p\in G: pe_0=\lambda e_0, \lambda >0\}.$$

\par We will now be considering jets (to infinite order), at $e_0$,
 of homogeneous
 functions and tensor fields. Homogeneity here refers to the usual
 dilations of the space $W=\mathbb{R}^{n+2}$.
We denote by $E (a)$ the space of jets
 of functions $u_a$, homogeneous of degree $a$ and by
$E^{IJ\dots M}(a)$ (there are $m$ indices $I,J,\dots ,L$) the
space of jets of functions of homogeneity $a$ that take values in
 the space $W\otimes\dots \otimes W$ (we are tensoring $m$ times).
For example $E^{IJK}$ is just a convenient way of recording that the
 jets take values in $W\otimes W\otimes W$. Moreover, ${E_I}^{JK}$
is the space of jets of functions that take values in
$W^{*}\otimes W\otimes W$. We denote by $H(a)$ the space of jets
at $e_0$ of homogeneous {\it harmonic} functions $u_a$. {\it
Harmonic} here means with respect
 to the operator $\Delta= \tilde{g}^{IJ}\partial^2_{IJ}$.

\par We define the space $K$ of
{\it jets of linearized ambient curvature
 tensors} around a point $e_0\in \mathbb{R}^{n+2}$  to be the set of
jets to infinite order of
4-tensor fields $\rho_{IJKL}$, with homogeneity $-2$, that satisfy the
 identities:

$$\rho_{[IJ]KL}=0, \text {} X^L\rho_{IJKL}=0\text{, }
\partial_{[H}\rho_{IJ]KL}=0,\text{} \rho_{[IJK]L}=0, \rho_{IJ[KL]}=0$$
and moreover each tensor $\partial^{(m)}_{A\dots C}\rho_{IJKL}$ is
trace-free with respect to the quadratic form $\tilde{g}$. Here
$[\dots ]$ stands for summation over all cyclic permutations of
the indices inside the brackets.

We then define
 the function $Eval$, that evaluates each such jet at $e_0$. Note that
$X=(X^I)$ is an element of $E^I(1)$.
We then denote $e^I=Eval(X^I)$. We recall
that coordinate differentiation defines a $P$-invariant map:

$$\partial_I: E^{JK\dots M}(s)\longrightarrow E_I^{JK\dots M}(s-1)$$

\par Now, let us define $lin\Pi (H(a),K)$ to be the polynomial in
 $H(a)\oplus K$
that arises from $\Pi (\tilde{u}_a,\tilde{R})$ by replacing each
 factor $\tilde{\nabla}^{(m)}_{AB\dots D}\tilde{R}_{IJKL}$ by a
 factor $\partial^{(m)}_{AB\dots D}\rho_{IJKL}$
and each factor $\tilde{\nabla}^r_{AB\dots D}\tilde{u}_a$ by a
factor $\partial^r_{AB\dots D}u_a, u_a\in H(a)$. For each $p\in P$, we
 define $plin\Pi (H(a),K)$
to stand for the complete contraction that arises from $lin\Pi
(H(a),K)$ in the following way: Let $p=(q^i_j)$ be
 in the form (\ref{transPhi}) and $\lambda$ be the top-left
component of the matrix $p$. Denote the index $\infty$ by $n+1$
 for convenience. Then $plin\Pi (H(a),K)$ arises from $lin\Pi (H(a),K)$
by substituting each factor $\partial^{(m)}_{AB\dots
D}\rho_{IJKL}$ by a factor $\lambda^2\sum_{A',B',\dots
,D',I',J',K',L'=0}^{n+1}\partial^{(m)}_{A'B'\dots
D'}\rho_{I'J'K'L'}q^{A'}_{A} \dots q^{L'}_{L}$. We also replace
each factor $\partial^r_{AB\dots D} u_a$ by a factor $\lambda^a
\sum_{A',B',\dots ,D'=0}^{n+1}\partial^r_{A'B'\dots D'}u_a
q^{A'}_{A}\dots q^{D'}_{D}$ and each factor $\tilde{g}^{IJ}$ by
$\lambda^{-2}\sum_{I',J'=0}^{n+1}\tilde{g}^{I'J'}q^I_{I'}q^J_{J'}$.

\par We recall the following fundamental fact, that follows
 from \cite{f:pitca}:

\begin{proposition}
\label{fefeprop} In the notation above, let us suppose that we can
show that if $lin\Pi (H(a),K)$ is $P$-invariant, $lin\Pi
(H(a),K):H(a)\oplus K\longrightarrow E(b)$,
 then $lin\Pi (H(a),K)$ can be written as:

$$lin\Pi (H(a),K)=\sum_{h\in H} a_h C^h(H(a), K)$$

where each $C^h(H(a), K)$ is a complete contraction (in
$\tilde{g}$), in the form:

\begin{equation}
\label{lincontr}
\begin{split}
&contr_{\tilde{g}}(\partial^{(m_1)}_{AB\dots
D}\rho_{IJKL}\otimes\dots\otimes
\partial^{(m_s)}_{A'B'\dots D'}\rho_{I'J'K'L'}\otimes
\\& \partial^{(p_1)}_{FG\dots H}u_a\otimes\dots\otimes
 \partial^{(p_o)}_{F'G'\dots H'}u_a)
\end{split}
\end{equation}
with $\sum_{i=1}^s( m_i+2)+\sum_{i=1}^o p_i=K$ and $o=q$.
\newline

\par It then follows that $\Pi(\tilde{R},\tilde{u}_a)$ can be
 written in the form:

$$\Pi(\tilde{R},\tilde{u}_a)=\sum_{h\in H'} a_h
\tilde{C}^h_{\tilde{g}} (\tilde{u}_a)$$ where each
$\tilde{C}^h_{\tilde{g}} (\tilde{u}_a)$ is in the form
(\ref{contraction}) with $\sum_{i=1}^s( m_i+2)+ \sum_{i=1}^o
p_i=K$ and $o=q$.
\end{proposition}

\par In view of the above, the rest of this section will focus on
 proving the hypothesis of the above Proposition. It will prove
useful to establish two isomorphisms between the space of jets at
$e_0$ of homogeneous harmonic functions and linearized ambient
curvatures, and  the space of two lists of tensors, which we will
denote by $H_{list}(a), K_{list}$. We use Propositions 1.2 and 4.1
from \cite{beg:itccg} to establish these two isomorphisms.

\par Proposition 1.2 in \cite{beg:itccg} implies that if $a\notin\mathbb{Z}_{+}$
(which is the case here, since we are assuming
$a=w\notin\mathbb{Z}_{+}$) the space $H(a)$ of jets at $e_0$ of
$a$-homogeneous harmonic functions, is
 isomorphic to the $P$-module  of lists:

\begin{equation}
\label{H} H_{list}(a)=\{ T^0,T^1,T^2,\dots \}, T^l\in
\odot^l_0W^{*}\otimes \sigma_{a-l}
\end{equation}
where $e^I(T^{l+1})_{IJ\dots M}=(a-l)(T^l)_{J\dots M}$. Here
$\sigma_q$ is the 1-dimensional representation of $P$ where the
element of $P$ in (\ref{transPhi}) is mapped to $\lambda^{-q}$.
 Also, $\odot^l$ stands for the symmetrized $l$-tensor product
 and $\odot^l_0$ stands for the trace-free part of the symmetrized
$l$-tensor product.
These conditions reflect the fact that we are dealing
with densities of weight $a$, the trace-free restriction reflects the harmonicity
of $u_a$ and the  restriction on the contraction against $e^I$
reflects the Euler homogeneity relations.

We also recall the Proposition 4.1 in \cite{beg:itccg} which shows
that that $K$ is isomorphic to a certain $P$-module of lists of
tensors, denote it by $K_{list}$, with symmetries and
anti-symmetries that model the usual properties of the curvature
and its covariant derivatives. We refer the reader to that paper.
To avoid confusion, we will denote the tensors $T^{(l)}$ in
Proposition 4.1 in \cite{beg:itccg} by $Q^l$.
\newline

\par Now, we will revert from thinking of the vector space of jets to thinking of
$P$-modules of lists of tensors. We thus study a polynomial
$\Pi(H_{list}(a),K_{list})$ in elements of the lists $H_{list}(a),
K_{list}$.

\par Now, by Weyl's classical invariant theory and using the fact that
$O(g(x_0))\subset P$, we have that $\Pi(H_{list}(a),K_{list})$ can
be written as:

$$\Pi(H_{list}(a),
K_{list})=\Pi_{even}(H_{list}(a),K_{list})$$
where $\Pi_{even}(H_{list}(a),K_{list})$ is a linear combination of
 complete contractions in the form:

\begin{equation}
\label{cpolu} contr_g(T^{k_1}\otimes\dots\otimes T^{k_s}\otimes
Q^{l_1} \otimes\dots\otimes Q^{l_r})
\end{equation}
the contractions are with respect to the metric $g$.

\par This reflects the fact that a conformally invariant differential
operator is still a Riemannian operator.
\newline

\par Now, recall from \cite{beg:itccg} the notion of a weak Weyl
invariant: Assume that $C$,

$$C: H_{list}(a)\oplus K_{list}\longrightarrow \odot^{(m)}_0 W\otimes \sigma_{b+m}$$
is a linear combination of partial contractions (with respect to
the metric $\tilde{g}$) of the tensors $Q^l$, $T^k$, $e^I$ (which
we denote by $e$, for short). If $C$ is of the form:

$$C=e\otimes\dots\otimes e\otimes   I$$
where the factor $e$ is tensored $m$ times and $I$ is a
polynomial,

$$I: H_{list}(a)\oplus K_{list}\longrightarrow \sigma_b$$
then we observe that $I$ must be $P$-invariant. (This follows because
$Pe=\lambda e$). We will then call $I$ a {\it weak Weyl invariant}.

\par So the $m$-tensor $C^{AB\dots D}$ has $C^{0\dots 0}=I$ and all the
 other components vanish. Equivalently, the tensor $C_{AB\dots}$ has
$C_{\infty\dots \infty}=I$ and all the other components vanish. We then claim:

\begin{proposition}
\label{RtowW}
Our polynomial $\Pi_{even}(H_{list}(a),K_{list})$ can be written as a
 weak Weyl invariant.
\end{proposition}

{\it Proof:} Using the quadratic form $\tilde{g}$:

$$\tilde{g}_{IJ}X^IX^J=g_{ij}X^iX^j + 2X^0X^\infty,$$
we see that $\Pi_{even}(H_{list}(a),K_{list})$ can be written as linear
 combination of
 complete contractions (in $\tilde{g}$) of the tensors
$T^l_{AB\dots D\infty\dots\infty}$, $Q^r_{IJKL,AB\dots
D\infty\dots\infty}$, $Q^r_{IJK\infty,AB\dots D\infty\dots\infty}$
and $Q^r_{I\infty K\infty,AB\dots D\infty\dots\infty}$ (where the
 indices $A,B,\dots ,D,I,J,K,L$ take the values $0,1,\dots ,n,\infty$).

By juxtaposing $e_I=(0,\dots ,0,1)$ if necessary, we may assume that
 the number of $\infty$ indices is equal to $m$ in all the terms in that
 linear combination.

\par Then, we make all the $\infty$'s into free indices $X,Y,\dots ,Z$
 and so we have an $m$-tensor $F_{XY\dots Z}$. We then take the symmetric
  and trace-free part of $F_{XY\dots Z}$, say $C_{XY\dots Z}$.
  Because $\tilde{g}_{\infty\infty}=0$, we still have that
$C_{\infty\dots\infty}=\Pi_{even}(H_{list}(a),K_{list})$.
Equivalently, raising indices we have that  $C^{0\dots
0}=\Pi_{even}(H_{list}(a), K_{list})$.

\par We then consider $D_{XY\dots Z}= C_{XY\dots Z}-
e\otimes\dots\otimes e\otimes \Pi_{even}(H_{list}(a),K_{list})$. All we
 need to show
is that $D_{XY\dots Z}$ vanishes. In other words, this is a linear algebra
 problem: Given the quadratic form $\tilde{g}$ and an $m$-form $C^{XY\dots Z}$ that is
 $P$-invariant and symmetric and totally trace-free so that
$C^{0\dots 0}=0$, then show that $C^{XY\dots Z}=0$. But this is shown
in Proposition 2.1 in \cite{beg:itccg} (we only need the even case
here). $\Box$
\newline

\par What remains to be done is to show that the term $I$ in
$\Pi_{even}(H_{list}(a),K_{list})=e\otimes\dots\otimes e\otimes I$
is a Weyl invariant. But this part of the argument exactly follows
the ideas in \cite{beg:itccg}:

\par We  now revert to thinking of the modules $H(a)$ and $K$ as
jets of homogeneous harmonic functions and linearized curvature tensors
around $e_0$, rather than just lists of tensors.

\par Let $F^{AB\dots G}(q+m)$ denote the space of jets at $e_0$ of
restrictions of $m$-tensors (with homogeneity $q+m$) to the light cone
 $Q$.

\par We now use the differential operator $D_I$ from \cite{beg:itccg}
that acts on elements of $F^{AB\dots E}(s)$ and maps them into
$F_I^{AB\dots E}(s-1)$. We pick any $f^{AB\dots E}\in F^{AB\dots E}(s)$
and  arbitrarily extend it
off of $Q$ to an element of $E^{AB\dots E}(s)$. We then define:

\begin{equation}
\label{formula}
D_If^{AB\dots E}=(\partial_I-\frac{X_I\Delta f^{AB\dots E}}
{(n+2s-2)})|_{Q}
\end{equation}

\par It follows that this operation is independent of the extension off of $Q$ (see \cite{beg:itccg}).
Moreover, we observe that for any $f\in F(s), n+2s\ne 0$, we have:

\begin{equation}
\label{grafe}
D_I(X^If)=\frac{(n+2s+2)(n+s)}{(n+2s)}f
\end{equation}

\par Repeatedly applying this operator as in \cite{beg:itccg}, we
then conclude that any weak Weyl invariant with $b=(\kappa\cdot
w-K)+m, m\in \mathbb{Z}_{+}$, with $w$ subject to the restrictions
of Theorem \ref{easy},
 will be a Weyl invariant. We have proven our Theorem
 \ref{easy}. (Note: The restriction
 $(w+\frac{n}{2})\notin\mathbb{Z}$ ensures that the operation is
 well-defined; The last two restrictions ensure that the constant
 on the right hand side of (\ref{grafe}) is not zero).
$\Box$

\section{Appendix: Proof of Lemma \ref{appendix}.}

 We first prove our claim for
complete contractions in the form (\ref{contr1}). Let us denote by
$\omega$ the number of factors of $C_g(f)$. It is straightforward
to observe that if $C_g(f)$ satisfies the extra restrictions then
all contractions $C^k_g(f)$ with $\omega$ factors (in the notation
of Lemma \ref{appendix}) will also satisfy the extra restrictions.
Thus, we may restrict attention to the contractions $C^k_g(f)$
with $\omega+1$ factors.

\par By definition, such contractions can arise in this setting
only by applying the curvature identity. But then our claim
follows because whenever we apply the curvature
 identity to a factor $F_h=\nabla^{(p)}f$ with $\beta[F_h]=\chi$,
or a factor $F_h=\nabla^{(m)}R$ with $\gamma[F_h]=\chi$, we will
 obtain a linear combination of partial contractions
in the forms $\nabla^{(m')}R\otimes\nabla^{(p')}f$ or
$\nabla^{(m')}R\otimes\nabla^{(m'')}R$, respectively. Denote these
factors by $F'_h, F''_h$. It follows that we will then have
$\beta[F''_h]\le \chi-2$ and $\gamma[F'_h]\le \chi$. Thus, since
we stated off with a complete contraction that satisfied the extra
restrictions, we  obtain a complete contraction that satisfies the
extra  restrictions.
\newline

\par Now we prove our claim for complete contractions in the form
(\ref{Wcontr}). In that setting, our claim is obvious if
$C^h_g(f)$ arises from $C_g(f)$ by switching two derivative
indices, or by applying (\ref{cotton}) to a factor $\nabla^{(p)}P$
or by applying a fake second Bianchi identity provided we {\it do
not} increase the number of factors and also {\it do not} increase
the number of internal contractions. The remaining cases in which
a contraction $C^k_g(f)$ may arise from $C_g(f)$ are as follows:

 Firstly, $C^k_g(f)$ may arise by applying
the Ricci identity to two derivative indices in the same factor.
Secondly  by replacing a factor $\nabla^{(m)}W$ by one of the
summands $\sum(\nabla^sW_{sjkl}\otimes g)$,
$\sum(\nabla^{ik}W_{ijkl}\otimes g)$ in the right hand sides of
(\ref{in&out}), (\ref{in&out2}). Thirdly, by the expressions
$Q(R)$ in the right hand sides of (\ref{in&out2}),
(\ref{in&out2c}).

\par Now, in the first case our claim follows by the same
argument as for $C_g(f)$ being in the form (\ref{contr1}). For the
second case we observe that whenever we apply (\ref{in&out}) or
(\ref{in&out2}) to a factor $F_h$ and we introduce a
 correction term of the form $F'_h=(\nabla^sW_{srty}\otimes g)$ or of the form
 $F'_h=(\nabla^{st}W_{srty}\otimes g)$ respectively,
  we obtain a complete contraction with at least $\mu+1$ internal
 contractions, and where the number of derivatives on each
 factor remains unaltered. Thus, we only have to check that the
 extra restrictions continue to hold whenever they are applicable.
  To see this, we observe by definition that if
 we had $\gamma[F_h]=\chi$ before we applied (\ref{in&out}) or
 (\ref{in&out2}), we then have $\gamma[F'_h]=\chi-1$. Furthermore,
 if we create an internal contraction in some other factor $T$ by
 virtue of the un-contracted metric tensor, we again decrease the
 quantities $\gamma[T]$, $\beta[T]$. All the other factors $T'$ will
 still have the same numbers $\beta[T'], \gamma[T']$, so we derive
 our claim.

\par The third case is when $C^h_g(f)$ arises from
$C_g(f)$ when we apply one of the equations
 (\ref{in&out2}), (\ref{in&out2c}) to a factor
$F_h=\nabla^{(m)}W$ and bring out the correction terms in $Q(R)$.
Suppose that $C^k_g(f)$ arises by applying (\ref{in&out2}) or
(\ref{in&out2c}) to a factor $F_h=\nabla^{(m)}W$.
 We only have to check that $F_h$ is then being replaced by a partial
contraction $\nabla^{(m')}R\otimes \nabla^{(m'')}R$ for which
$\gamma[\nabla^{(m')}R]<\frac{n}{2}$,
$\gamma[\nabla^{(m')}R]<\frac{n}{2}$, whenever this extra
restriction is applicable. Now, for this we see that $Q(R)$ is a
sum of partial contractions in the form $R\otimes R$. But each
such curvature expression has $\gamma[R]\le 2$, whereas the
expression $F_h$ in the left hand sides of (\ref{in&out2}),
(\ref{in&out2c}) have $\gamma[F_h]=3$. Furthermore, we observe
that $m'+m''=m-2$ and that each derivative index $r$ in $F_h$ that
does not belong to an internal contraction in $F_h$ will belong to
only one of the factors $\nabla^{(m')}R, \nabla^{(m'')}R$, while
each internal contraction in $F_h$ will either give rise to an
internal contraction in one of the factors
$\nabla^{(m')}R,\nabla^{(m'')}R$, or it will give rise to two
indices ${}_a,{}_b$ in $\nabla^{(m')}R,\nabla^{(m'')}R$
respectively that will contract against each other. Thus, since
$F_h$ satisfied $\gamma[F_h]<\frac{n}{2}$, we obtain that
$\gamma[\nabla^{(m')}R]<\frac{n}{2}$,
$\gamma[\nabla^{(m')}R]<\frac{n}{2}$. $\Box$


\begin{thebibliography}{12}

\bibitem{a:ocidood} S. Alexakis \emph{On conformally invariant differential operators in odd dimensions}
 Proc. Natl. Acad. Sci. USA  {\bf 100}  (2003),  no. 8, 4409--4410


\bibitem{a:dgciI} S. Alexakis \emph{The decomposition of Global
Conformal Invariants I} to appear in Ann. of Math.

\bibitem{a:thesis} S. Alexakis \emph{PhD Thesis}, Princeton University 2005.


\bibitem{b:tdcpitci} T. N. Bailey \emph{Thomas' $D$-calculus,
Parabolic Invariant theory, and Conformal Invariants}, Twistor
Theory (Plymouth), 1-7, Lecture Notes in Pure and Appl. Math. {\bf169}, Dekker, New York, 1995.

\bibitem{beg:tsbcprs} T.N. Bailey, M. G. Eastwood, A. R. Gover
\emph{Thomas' structure bundle for conformal, projective and
related structures} Rocky Mountain J. Math {\bf 24} (1994),
1191-1217.

\bibitem{beg:itccg} T. N. Bailey, M. G. Eastwood, C. R. Graham
\emph{Invariant Theory for Conformal and CR Geometry} Ann. of Math (2),
 {\bf 139} (1994), 491-552.


\bibitem{b:fd} T. Branson \emph{The functional determinant}, Global
Analysis Research Center Lecture Note Series, no. 4, Seoul
National University (1993).





\bibitem{eg:icd} M.G. Eastwood, C.R. Graham \emph{Invariants of
Conformal Densities},  Duke Math. J.  {\bf 63}  (1991),  no. 3,
633-671.


\bibitem{e:ntrm} D.B.A. Epstein \emph{Natural Tensors on Riemannian
Manifolds}, Journal of Diff. Geom. {\bf 10} (1975), 631-645


\bibitem{f:pitca} C. Fefferman \emph{Parabolic Invariant Theory in
 Complex Analysis.}, Adv. in Math. {\bf 31} (1979), 131-262.

\bibitem{fg:ci} C. Fefferman, C. R. Graham \emph{Conformal Invariants}
\'Elie Cartan et les mathematiques d'aujourd'hui,      Ast\'erisque,
numero hors serie, 1985, 95-116.

\bibitem{fh:amcqcccg} C. Fefferman, K. Hirachi \emph{Ambient Metric
Construction of Q-Curvature in Conformal and CR Geometries}, Math. Res.
Lett. {\bf 10} (2003), 819-831.


\bibitem{g:itpg} C. R. Graham \emph{Invariant theory of Parabolic
Geometries},  Complex geometry (Osaka, 1990),  53--66,
Lecture Notes in Pure and Appl. Math., 143, Dekker, New York, 1993.





\bibitem{gjms:cipl} C. R. Graham, R. Jenne, L. J. Mason, G. Sparling
\emph{Conformally invariant powers of the Laplacian. I. Existence},
J. London Math. Soc. (2) {\bf 46} (1992), 557-565.

\bibitem{g:ciplne}
C. R. Graham, \emph{Conformally invariant powers of the Laplacian. II. Nonexistence},
  J. London Math. Soc. (2)  {\bf 46}  (1992),  566--576.

\bibitem{g:itccg} A. R. Gover \emph{Invariant Theory and Calculus
for Conformal geometries}, Adv. Math. {\bf 163} (2001), 206-257.


\bibitem{gh:cipl} A.R. Gover, K. Hirachi \emph{Conformally
invariant powers of the Laplacian: A complete non-existence
theorem}, J. Am. Math. Soc. {\bf 17} (2004), 389-405.

\bibitem{h:new} K. Hirachi \emph{Ambient metric construction of CR invariant differential
operators}, preprint  math.CV/0701804.


\bibitem{j:cc} H. P. Jakobsen \emph{Conformal covariants},
 Publ. Res. Inst. Math. Sci. 22 (1986), no. 2, 345–364.

\bibitem{k:tgdg} S.  Kobayashi, \emph{Transformation groups in differential geometry},
Springer-Verlag (1972).

\bibitem{p:cti} R. Penrose, \emph{A Conformal treatment of Infinity}, 1964
 Relativit\'e, Groupes et Topologie (Lectures, Les Houches, 1963 Summer School of Theoret. Phys., Univ. Grenoble)

\bibitem{s:ct} P. Szekeres \emph{Conformal tensors}, Proc. Roy.
Soc. London Ser. A., {\bf 304}, 1968, 113-122.


\bibitem{t:cg} T. Y. Thomas, \emph{On Conformal Geometry}, Proc. Nat. Acad. Sci. USA {\bf 12}, (1926),
352-359


\bibitem{t:digs} T. Y. Thomas, \emph{The differential invariants of generalized spaces},
Cambridge University Press, (1934)










\bibitem{w:cg} H. Weyl \emph{The Classical Groups}, Princeton University Press 1946.


\end{thebibliography}
\end{document}